\documentclass[preprint]{elsarticle}
\usepackage{amsfonts,amssymb}        %
\usepackage{amsmath}
\usepackage{latexsym}
\usepackage{verbatim}    %
  
\usepackage{bm,enumerate}

\usepackage{multirow}
\usepackage{array}

\usepackage{float}
\usepackage{graphics}
\usepackage{subfigure}
\usepackage{color}
\usepackage{wrapfig}
\usepackage{sidecap}

\def\SMALLO{\mathcal{o}}
\def\DELTAT{{\delta t}}
\def\NUMSP{{K}}                       

\def\INDXE{\mathcal{I}^E}
\def\INDXO{\mathcal{I}^O}

\def\SINH{\mathrm{sinh}}
\def\COSH{\mathrm{cosh}}

\def\BNDRY{\partial}
\def\CLSR{\bar}

\def\DIAM{\mathrm{diam}\,}

\def\RANG_S{S}                            

\def\LATT{{\Lambda}_N}                    
\def\SIGMA{{\mathcal{S}}}               

\def\SPINSP{\Sigma}                     
                   
\def\CONFNEW{\sigma^{x,\omega}}

\def\PROCMIC{\{S_t\}_{t\geq 0}}

\def\PROC#1{\{{#1}\}_{t\geq 0}}

\def\CUBE{C}

\def\LOPER{\mathcal{L}}

\def\Oo{\mathcal{O}}

\def\EXP#1{e^{#1}}

\def\EXPECT{{\mathbb{E}}}

\def\PROB#1{{\mathbb{P}\left({#1}\right)}}

\def\COMMA{\,,}             
\def\PERIOD{\,.}            
\def\SEP{{\,|\,}}           

\def\VIZ#1{(\ref{#1})}      

\def\BIGO{\Oo}
\def\SMALLO{\mathrm{o}}

\def\LTWO{L^2}
\def\CB{C_b}

\journal{ Journal of Computational Physics}

\begin{document}

\begin{frontmatter}

\title{Hierarchical fractional-step approximations and parallel kinetic Monte Carlo
algorithms}
\author[crete]{Giorgos Arampatzis}%
\ead{garab@math.uoc.gr}
\author[umass]{Markos A. Katsoulakis\corref{cor1}}%
\ead{markos@math.umass.edu}
\author[udel]{Petr Plech\'a\v{c}}%
\ead{plechac@math.udel.edu}
\author[udel-cs]{Michela Taufer}
\ead{taufer@cis.udel.edu}
\author[udel-cs]{Lifan Xu}
\ead{xulifan@udel.edu}

\address[crete]{Department of Applied Mathematics, 
                University of Crete and  Foundation of Research and Technology-Hellas, Greece }
\address[umass]{Department of Mathematics and Statistics, University of Massachusetts at Amherst, 
                Amherst, MA 01003, USA}
\address[udel]{Department of Mathematical Sciences,
               University of Delaware, Newark, DE 19716, USA }
\address[udel]{Department of Computer Science,
               University of Delaware, Newark, DE 19716, USA }

\cortext[cor1]{Corresponding author}

\begin{abstract} 
We present a mathematical framework  for constructing and analyzing parallel algorithms for lattice Kinetic Monte Carlo (KMC) simulations.
The resulting   algorithms have the capacity to simulate a wide range of spatio-temporal scales in  spatially distributed, non-equilibrium physiochemical 
processes with complex chemistry and transport micro-mechanisms. The algorithms can be tailored to specific hierarchical parallel 
architectures such as multi-core processors or clusters of Graphical Processing Units (GPUs). 
The proposed parallel algorithms are controlled-error approximations of kinetic Monte Carlo algorithms, departing from the predominant paradigm of 
creating  parallel KMC algorithms with exactly the same master equation as the serial one. 

Our methodology relies on %
a spatial decomposition of  the Markov operator underlying  the KMC algorithm  into  a hierarchy of operators 
corresponding  to the  processors' structure  in the parallel architecture.  Based on this operator decomposition, we  formulate  
{\it Fractional Step Approximation schemes}
by employing the Trotter Theorem and its random variants; these schemes,   
(a)  determine  the {\em communication schedule} between processors, and (b) are  run independently on each processor  through a serial KMC simulation, called a {\em kernel\,},
on each fractional step time-window. 
Furthermore,  the proposed mathematical framework allows us to rigorously justify  
 the numerical and statistical consistency of the proposed algorithms,  showing the convergence of our approximating 
schemes to  the original serial KMC. The approach also provides a systematic  evaluation of different processor 
communicating schedules.%
We carry out a detailed benchmarking of the parallel KMC schemes   using  available exact solutions, for example, in Ising-type systems and we demonstrate the 
capabilities of the method to simulate complex spatially distributed reactions at very large scales on  GPUs. 
Finally, we  discuss {\em work load balancing}  between processors and propose a re-balancing scheme based on probabilistic mass transport methods.

\end{abstract}


\end{frontmatter}
\section{Introduction}

Kinetic Monte Carlo algorithms have proved to be an important tool for the simulation of out-of-equilibrium, spatially distributed processes. 
Such models arise in physiochemical applications ranging from materials science and catalysis, to complex biological processes. 
Typically the simulated models involve chemistry and/or transport micro-mechanisms for atoms and molecules, e.g.,
reactions, adsorption, desorption processes and diffusion on surfaces and through  complex media,  
\cite{binder, Auerbach, Vlachos}. Furthermore, mathematically similar mechanisms and 
corresponding Kinetic Monte Carlo simulations arise in agent-based, evolutionary games problems in epidemiology, ecology and traffic networks, \cite{Szabo}.

The simulation of stochastic lattice systems using Kinetic Monte Carlo (KMC) methods relies on the direct numerical simulation of the underlying 
Continuous Time Markov Chain (CTMC). Since such stochastic processes are set on a lattice (square, hexagonal, etc.) $\LATT$ with $N$ sites, they have a 
discrete, albeit high-dimensional, configuration space  $\Sigma$ and necessarily have to be of jump type describing transitions between different 
configurations $\sigma \in \Sigma$. 
Mathematically, CTMC  are defined in terms of  the transition rates $c(x,\omega;\sigma)$ 
which correspond to an updating micro-mechanism that describes completely the evolution of the stochastic process as a transition from  
a current configuration $\sigma$ of the system to a new configuration $\CONFNEW$ by performing an update
in a neighborhood of  the site $x\in\LATT$.  
In other words the probability of a transition over an infinitesimal time interval $\DELTAT$ is
$\PROB{S_{t_\DELTAT} = \CONFNEW\SEP S_t = \sigma} = c(x,\omega;\sigma)\DELTAT + \SMALLO(\DELTAT^2)$.
In turn, the transition rates define the total rate
\begin{equation}\label{totalrate}
\lambda(\sigma)=\sum_{x \in \LATT}\sum_{\omega \in \SIGMA_x } c(x, \omega; \sigma)\COMMA
\end{equation}
which is the intensity of the exponential waiting time for a jump to be performed when the system is currently at the state $\sigma$. 
Here $\omega\in\SIGMA_x$, where $\SIGMA_x$ is the set of all possible
configurations that correspond to an update at a neighborhood $\Omega_x$ of the site $x$. 
Once this exponential ``clock'' signals a jump, then the system transitions from the state $\sigma$ to a new configuration $\CONFNEW$ with probability
\begin{equation}\label{skeleton}
p(\sigma, \CONFNEW)={c(x, \omega; \sigma)\over \lambda(\sigma)}\, .
\end{equation}
Thus the full stochastic evolution is completely defined. We refer to the discussion in Section~\ref{formulation} for a complete mathematical 
description of the KMC method.
The implementation of  this method is based on efficient calculation of \VIZ{totalrate} and \VIZ{skeleton}, and was first developed in \cite{BKL75}, known as a BKL Algorithm,  for 
stochastic lattice Ising models, and in \cite{Gillespie76} known as Stochastic Simulation Algorithm (SSA) for reaction systems.
However, as it is evident from formulas \VIZ{totalrate} and \VIZ{skeleton}, the algorithms are inherently serial as updates are done 
at one site $x \in \LATT$ at a time, while on the  other hand  the calculation of \VIZ{totalrate}  depends on information from the entire spatial 
domain $\LATT$. For these reasons it seems, at first glance, that KMC algorithms cannot be parallelized easily.

However, Lubachevsky, in \cite{Lubachevsky88}, proposed  an {\em asynchronous}  approach for parallel KMC simulation  in the context of Ising systems, 
in the sense that different processors simulate independently parts of the physical domain, while  inconsistencies 
at the boundaries are corrected with a series of suitable rollbacks. This method relies on uniformization of the total  rates over each processor, 
see also \cite{Nicol} for the use of uniformization in the parallel simulation of general CTMC. Thus the approach yields a {\em null-event} algorithm, \cite{binder},  which includes 
rejected moves over the entire domain of each processor.
Furthermore,  Lubachevsky proposed a modification in order to incorporate the BKL Algorithm in his parallelization method, which was implemented 
and tested in \cite{Korniss99}.  This is a partially rejection-free  (still asynchronous) algorithm, where  BKL-type  rejection-free simulations are carried out 
in the interior of each processor, while uniform rates were used at the boundary, reducing rejections over just the boundary set. However, in spite of the 
proposed improvements, these asynchronous algorithms may still have a high number of rejections for boundary events and rollbacks, which considerably  
reduce the parallel efficiency, \cite{ShimAmar05}. 
Advancing processors in time in a synchronous manner over a fixed time-window can provide a way to mitigate the excessive number of boundary inconsistencies 
between processors and ensuing rejections and rollbacks in earlier methods. Such {\em synchronous} parallel KMC algorithms were proposed and extensively 
studied in \cite{Lubachevsky93, ShimAmar05, MerickFichthorn07, ShimAmar09}. However, several costly global communications are required at each cycle 
between all processors, whenever a boundary event occurs in any one of them, in order to avoid errors in the inter-processor communication and rollbacks, \cite{ShimAmar09}.

As we will discuss further in this paper, many of the challenges in parallel KMC  can be addressed by abandoning the 
earlier perspective on creating a parallel KMC algorithm with the exactly same rates (and hence the generator and master equation) as the serial algorithm, 
see \cite{Kalos} for a discussion on exact algorithms. This is a very natural idea in the numerical analysis of continuum models such as 
Ordinary and Partial Differential Equations (ODE/PDE). First, in \cite{ShimAmar05b} the authors propose an {\em approximate} algorithm,   
in order to create a parallelization scheme for KMC. It was recently demonstrated \cite{ShimAmar09, SPPARKS}, that this method is very promising:  
boundary inconsistencies  are resolved in a straightforward fashion, while there is an absence of global communications in contrast to synchronous 
relaxation schemes discussed earlier. Finally, we note that, among the parallel algorithms  tested in \cite{ShimAmar09}, the  approximate  algorithm had  
the highest parallel efficiency. 

Here we  develop a general mathematical framework for {\em parallelizable  approximations} of the KMC algorithm.%
Our approach  relies on first  developing a spatial decomposition of the Markov operator, that defines  the Kinetic 
Monte Carlo algorithm,  into  a hierarchy of operators. The decomposition is tailored   to the  processor architecture. 
Based on this operator 
decomposition, we  formulate {\it Fractional Step Approximation} schemes
by employing the Trotter product formula. In turn these approximating schemes  determine {\em Communication Schedule} between processors through 
the sequential application of the operators in the decomposition, and the time step employed in the particular fractional 
step scheme.
Here we discuss deterministic schedules resulting from Lie- and  Strang-type fractional step schemes, 
as well as random schedules derived 
by  the Random Trotter Theorem, \cite{Kurtz}. 
We show that the scheme in \cite{ShimAmar05b} is a particular case of a random schedule
and can be mathematically  analyzed within the proposed framework.%
We recall that the deterministic Trotter Theorem was first proved in \cite{Trotter} for the approximation 
of semigroups corresponding to operator sums, and it has found wide application in the numerical ODE/PDE analysis, e.g., 
\cite{Hairer}.

In Section~\ref{formulation} we show that the Fractional Step KMC  schemes   allow  us to run independently on each 
processor  a serial KMC simulation (called a {\em kernel\,})
on each fractional time-step window. Furthermore, processor communication is straightforward at the end of each 
fractional time-step while  no global communications or rollbacks  are  involved. 
In Section~\ref{gpu} we show that the hierarchical structure of our methodology can be easily  implemented for 
very general physiochemical processes modeled by lattice systems,  allowing    users to  input  as  the algorithm's
KMC kernel their preferred serial algorithm. This flexibility and hierarchical structure are
 key  advantages
for tailoring our framework to particular parallel architectures   with complex  memory and processor hierarchies, 
e.g., clusters of GPUs.

The proposed mathematical framework allows  us to rigorously prove   
the numerical and statistical consistency of the proposed algorithms,   while on the other hand it 
provides a systematic  evaluation of different processor communication schedules.
Indeed, in Section~\ref{pcs}  the numerical and statistical consistency of the proposed algorithms is rigorously 
justified by the Trotter Theorem, \cite{Trotter}, \cite{Hairer} showing the convergence 
of our approximating schemes to the original serial KMC algorithm, interpreted as convergence to the underlying  
Markov operator. 
Using the Random Trotter Theorem  \cite{Kurtz} we show that
the approximation schemes with a random schedule, including the one in  \cite{ShimAmar05b} as a special case, 
are  numerically consistent in the approximation limit; that is,  as the time step in the fractional step scheme 
converges to zero, it converges to a continuous time Markov Chain that has the same master equation and generator 
as the original serial KMC. 
 In Section~\ref{cerror}  we show that the proposed mathematical framework can allow the  study of 
controlled-error approximation properties of Fractional Step KMC schemes, as well as the systematic  evaluation 
of different processor communicating schedules, comparing  for instance the scheme in \cite{ShimAmar05b} to the 
Lie scheme \VIZ{lie}.
Finally, in Section~\ref{workloadb} we  discuss work-load balancing between processors and propose a re-balancing scheme 
based on probabilistic mass transport methods,
\cite{lcevans}, which is particularly well-suited for the proposed fractional step KMC methods.
In Section~\ref{apps}  we present  detailed benchmarking of the proposed parallel algorithms  using analytically available 
exact solutions, for instance, in Ising-type systems and demonstrate the 
capabilities of the method to simulate  complex spatially distributed molecular systems, 
such as CO oxidation on a catalytic surface.
\section{Fractional Step Kinetic Monte Carlo Algorithms}\label{formulation}
We first present the mathematical background of KMC in a more abstract way in order to demonstrate the 
generality and the applicability of the proposed method. We consider a $d$-dimensional lattice $\LATT$ with $N$ lattice sites.
We restrict our discussion to  lattice gas models where the order
parameter or the spin variable takes value in a finite countable set $\SPINSP=\{0,1,\dots, \NUMSP \}$.
At each lattice site $x\in \LATT$ an order parameter (a spin variable) $\sigma(x)\in \SPINSP$ 
is defined. The states in $\SPINSP$ correspond to occupation of the site $x\in\LATT$ by different
species. For example, if $\SPINSP=\{0,1\}$ the order parameter models the classical lattice gas with
a single species occupying the site $x$ when $\sigma(x)=1$ and with the site being vacant if $\sigma(x)=0$.
We denote $\PROCMIC$ the stochastic process with values in the configuration space
 $\SIGMA=\SPINSP^{\LATT}$.

Our primary focus is on modeling the basic processes of adsorption, desorption, diffusion and reactions
between different species. Thus the local dynamics is described by a collection of the transition rates
$c(x,\omega;\sigma)$ and by an updating mechanism such that 
the configuration $\sigma$ of the system changes into a new configuration $\CONFNEW$ by an update
in a neighborhood of  the site $x\in\LATT$. Here $\omega\in\SIGMA_x$, where $\SIGMA_x$ is the set of all possible
configurations that correspond to an update at a neighborhood $\Omega_x$ of the site $x$. For example, 
if the modeled process is 
a diffusion of the classical lattice gas a particle at $x$, i.e., $\sigma(x)$ can move to any 
nearest neighbor of $x$, i.e., $\Omega_x = \{y\in\LATT\SEP |x-y|=1\}$ and $\SIGMA_x$ is the set of all possible
configurations $\SIGMA_x = \SPINSP^{\Omega_x}$. In other words the collection of
measures $c(x,\omega;\sigma)$ defines the transition probability  from $\sigma$ to $\CONFNEW$ over an infinitesimal
time interval $\DELTAT$. More precisely, the evolution of the system is described by a continuous time Markov jump process with the
generator $\LOPER:\CB(\SIGMA) \to \CB(\SIGMA)$ acting on  continuous bounded test functions $f \in \CB(\SIGMA)$ according to
\begin{equation}\label{generator}
 \LOPER f(\sigma) = \sum_{x\in\LATT}\sum_{\omega\in\SIGMA_x} c(x,\omega;\sigma) [f(\CONFNEW) - f(\sigma)]\PERIOD
\end{equation}
We recall that the evolution
of the expected value for an arbitrary observable $f \in \CB(\SIGMA)$ is given by the action of the Markov semigroup $\EXP{t\LOPER}$ associated with the generator $\LOPER$ and the process $\PROC{S_t}$
\begin{equation}\label{semigroup}
 \langle \EXP{t\LOPER}\mu_0, f\rangle = \EXPECT_{S_0}[f(S_t)]\COMMA
\end{equation}
where $\mu_0$ is the initial distribution of the process, i.e. of the random variable $S_{0}$, \cite{Liggett}.
Practically, the sample paths $\PROCMIC$
are constructed 
via KMC, that is  through the procedure described in \VIZ{totalrate} and \VIZ{skeleton}. 

To elucidate the introduced notation we give a few examples relevant 
to the processes modeled here. We refer, for instance,  to \cite{binder, Auerbach, Vlachos} for a complete discussion of the physical processes.

\medskip

\noindent{\sc Examples.}
\begin{enumerate}
 \item {\it Adsorption/Desorption for single species particles.} In this case spins take values in
       $\sigma(x)\in\SPINSP=\{0,1\}$, $\Omega_x = \{x\}$, $\SIGMA_x = \{0,1\}$ and the update represents a spin flip
       at the site $x$, i.e., for $z\in\LATT$
       $$
	  \CONFNEW(z) \equiv \sigma^x(z) =  \begin{cases} 
	                                       \sigma(z) & \mbox{if $z\neq x$,} \\
                                               1 - \sigma(x) & \mbox{if $z = x$.}
	                                    \end{cases}
       $$
 \item {\it Diffusion for single species particles.} The state space for spins is $\sigma(x)\in\SPINSP=\{0,1\}$, 
       $\Omega_x = \{y\in\LATT\SEP |x-y|=1\}$ includes all nearest neighbors of the site $x$ to which a particle
       can move. Thus the new configuration $\CONFNEW = \sigma^{(x,y)}$ is obtained by updating the configuration $S_t=\sigma$ 
       from the set of possible local configuration changes $\{0,1\}^{\Omega_x}$ using the specific rule, also known
       as spin exchange, which involves changes at two sites $x$ and $y\in\Omega_x$
       $$
	  \CONFNEW(z) \equiv \sigma^{(x,y)}(z) =  \begin{cases} 
	                                              \sigma(z) & \mbox{if $z\neq x,y$,} \\
                                                      \sigma(x) & \mbox{if $z = y$,} \\
                                                      \sigma(y) & \mbox{if $z = x$.}
	                                         \end{cases}
       $$
       The transition rate is then written as $c(x,\omega;\sigma) = c(x,y;\sigma)$.
       The resulting process $\PROCMIC$ defines  dynamics with the total number of particles ($\sum_{x\in\LATT}\sigma(x)$)
       conserved, sometimes referred to as Kawasaki dynamics.
\item {\it Multicomponent reactions.} Reactions that involves $\NUMSP$ species of particles are easily described 
      by enlarging the spin space to $\SPINSP=\{0,1,\dots,\NUMSP\}$. If the reactions occur only at a single site
      $x$, the local configuration space $\SIGMA_x = \SPINSP$ and the update is indexed by $k\in\SPINSP$ with the
      rule
      $$
	  \CONFNEW(z) \equiv \sigma^{(x,k)}(z) =  \begin{cases} 
	                                              \sigma(z) & \mbox{if $z\neq x,y$,} \\
                                                      k  & \mbox{if $z = x$.}
	                                         \end{cases}
      $$
      The rates $c(x,\omega;\sigma) \equiv c(x,k;\sigma)$ define probability of a transition $\sigma(x)$ to species
      $k=1,\dots,\NUMSP$ or vacating a site, i.e., $k=0$, over $\DELTAT$.

      \medskip

\item {\it Reactions involving  particles with internal degrees of freedom.}
      Typically a reaction involves particles with internal degrees of freedom, and in this case  several neighboring lattice sites may be updated at the same time, corresponding to the degrees of freedom of the  particles involved in the reaction. For example, in a case such as  CO oxidation on a catalytic surface, \cite{evans09}, 
      when only particles at a nearest-neighbor distance can react we set $\sigma(x)\in\SPINSP=\{0,1,\dots,\NUMSP\}$,
      $\Omega_x = \{y\in\LATT\SEP |x-y|=1\}$ and the set of local updates $\SIGMA_x = \SPINSP^{\Omega_x}$. Such $\SIGMA_x$
      contains all possible reactions in a neighborhood of $x$. When reactions involve only pairs of species, the rates can be
      indexed by $k$, $l\in \SPINSP$, or equivalently $\SIGMA_x = \SPINSP\times \SPINSP$. 
      Then the reaction rate $c(x,\omega;\sigma) = c(x,y,k,l;\sigma)$ describes the  probability per unit
      time of  $\sigma(x)\to k$ at the site $x$ and $\sigma(y)\to l$ at $y$, i.e., the updating mechanism
      $$
	  \CONFNEW(z) \equiv \sigma^{(x,y,k,l)}(z) =  \begin{cases} 
	                                              \sigma(z) & \mbox{if $z\neq x,y$,} \\
                                                      k         & \mbox{if $z = x$,} \\
                                                      l         & \mbox{if $z = y$,}
	                                         \end{cases}
      $$
      where $|x-y|=1$.
\end{enumerate}

\subsection{Hierarchical structure of the generator}
\begin{figure}
 \centerline{%
        \subfigure[]{\label{latt1}\includegraphics[angle=90,scale=0.55]{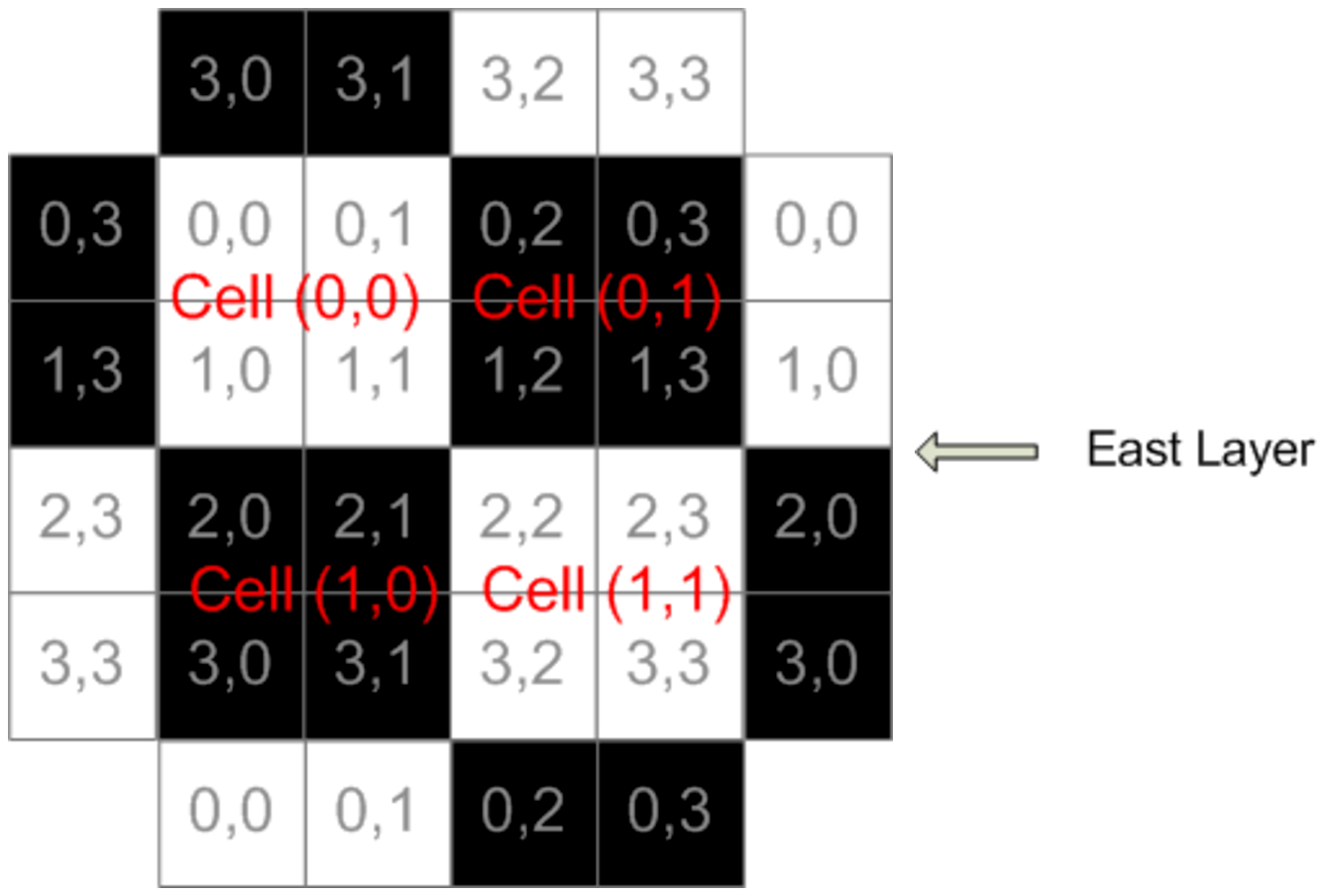}}%
        \hspace{-9cm}
        \subfigure[]{\label{latt2}\includegraphics[angle=90,scale=0.40]{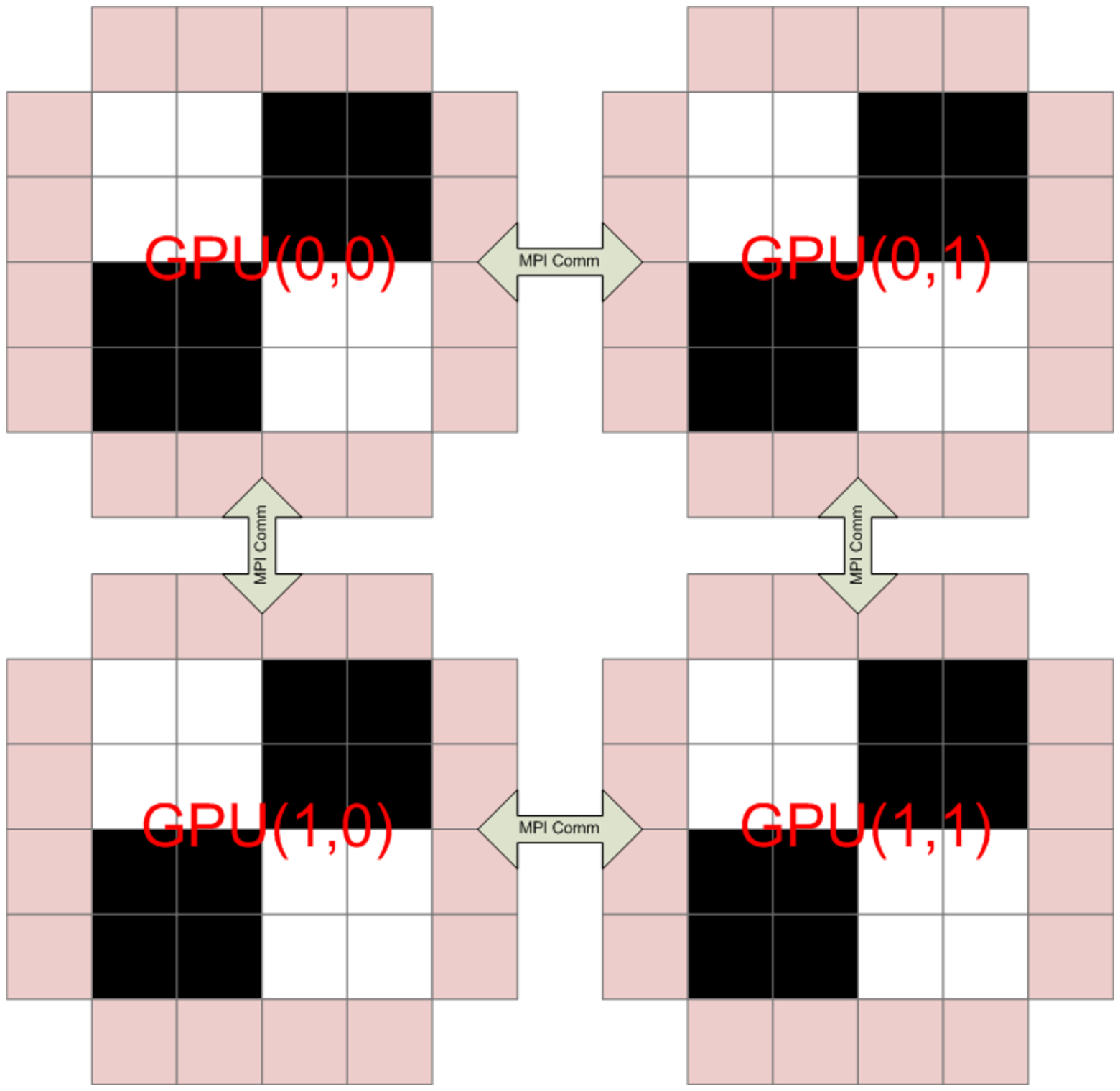}}    
                   }%
  \caption{(a) Lattice decomposition in \VIZ{sublatt} using the checkerboard scheme mapped onto a single
           multi-threading processing unit (e.g., GPU). 
           The integer cell coordinates also indicate communication through boundary buffer regions. In practice 
           other partitionings may result in a lower communication overhead.  (b) Hierarchical lattice partitioning on a cluster of processing units. 
  }
  \label{latt}
\end{figure}

The generator %
of the Markov process $\PROCMIC$ given in a general form in \VIZ{generator} is our starting point for the development 
of parallel algorithms based on geometric partitioning of the lattice. 
The lattice $\LATT$ is decomposed into non-overlapping cells $\CUBE_m$, $m=1,\dots,M$ such that 
\begin{equation}\label{decomposition}
 \LATT = \bigcup_{m=1}^M \CUBE_m\COMMA\;\;\;\; \CUBE_m \cap \CUBE_n = \emptyset\COMMA\; m\neq n\PERIOD
\end{equation}
With each set $\CUBE_m$ a larger set $\CLSR\CUBE_m$ is associated by adding sites to $\CUBE_m$ which are connected
with sites in $\CUBE_m$ by interactions or the updating mechanism,  see Figure~\ref{latt1}. More precisely, we define the range of
interactions $L$ for the set $\CUBE_m$ and the closure of this set
$$
\CLSR\CUBE_m = \{ z\in\LATT\SEP |z - x|\leq L\COMMA x\in\CUBE_m\}\COMMA\;\;
\mbox{where}\;L = \max_{x\in\CUBE_m}\{\DIAM \Omega_x\}\PERIOD
$$
In many models the value of $L$ is independent of $x$ due to translational invariance of the model.
The boundary of $\CUBE_m$ is then defined as
$\BNDRY\CUBE_m = \CLSR\CUBE_m\cap\CUBE_m$.  
This geometric partitioning induces a decomposition of the generator \VIZ{generator}
\begin{align}\label{gendecomp}
 \LOPER f(\sigma) & = \sum_{x\in\LATT}\sum_{\omega\in\SIGMA_x} c(x,\omega;\sigma) [f(\CONFNEW) - f(\sigma)] \\
                  & = \sum_{m=1}^M \sum_{x\in\CUBE_m} \sum_{\omega\in\SIGMA_x} c(x,\omega;\sigma) [f(\CONFNEW) - f(\sigma)]\\
                  & = \sum_{m=1}^M \LOPER_m f(\sigma)\PERIOD
\end{align}
The generators $\LOPER_m$ define new Markov processes $\{S^m_t\}_{t\geq 0}$ on the {\it entire} lattice $\LATT$.

 \medskip
 \noindent
 {\bf Remark:} In many models the interactions between particles are of the two-body type with the nearest-neighbor range
 and therefore the transition rates $c(x,\omega;\sigma)$ depend on the configuration $\sigma$ only through
 $\sigma(x)$ and $\sigma(y)$ with $|x-y|=1$. Similarly the new configuration $\CONFNEW$ involve changes only at the
 sites in this neighborhood. Thus the generator $\LOPER_m$ updates the lattice sites at most in the set 
 $\CLSR \CUBE_m = \{ z\SEP |x-z|=1\COMMA x\in\CUBE_m\}$, see Figure~\ref{latt1}. Consequently the processes $\PROC{S^m_t}$ and $\PROC{S^{m'}_t}$
 corresponding to $\LOPER_m$ and $\LOPER_{m'}$ are independent provided $\CLSR\CUBE_m \cap \CLSR\CUBE_{m'} = \emptyset$.

Therefore, splitting \VIZ{gendecomp} allows us to define independent processes which yields an algorithm suitable
for parallel implementation, in particular, in the case of short-range interactions when the communication
overhead can be handled efficiently. If the lattice $\LATT$ is partitioned into subsets $\CUBE_m$ such that
the diameter $\DIAM\CUBE_m > L$, where $L$ is the range of interactions, we can group the sets
$\{\CUBE_m\}_{m=1}^M$ in such a way that there is no interaction between sites in the sets $\CUBE_m$ that
belong to the same group.
For the sake of simplicity we assume that the lattice is divided into two {\em sub-lattices}   described
by the index sets $\INDXE$ and $\INDXO$, (black vs. white in Fig.~\ref{latt1}), hence we have
\begin{equation}\label{sublatt}
\LATT = \LATT^E\cup\LATT^O:=\bigcup_{m\in\INDXE} \CUBE_m^E \cup \bigcup_{m\in\INDXO} \CUBE_m^O\PERIOD
\end{equation}
Other lattice partitionings are also possible and may be more suitable for specific   
micro-mechanisms in the KMC or the computer architecture. Choice of the partitioning scheme can reduce communication 
overhead, see for instance \cite{ShimAmar05b}. 
For the sake of simplicity in the presentation, here we consider the partitioning depicted in \VIZ{sublatt} and 
Fig.~\ref{latt1}, although our mathematical framework applies to any other sublattice decomposition. %
Returning to \VIZ{sublatt}, the sub-lattices {\em induce}  a corresponding splitting of the generator:
\begin{equation}\label{opdecomp}
\LOPER =\LOPER^E+\LOPER^O:= \sum_{m\in\INDXE} \LOPER_m^E + \sum_{m\in\INDXO} \LOPER_m^O\PERIOD
\end{equation}
This simple observation has key consequences for simulating the process $\PROC{S_t}$ in parallel, 
as well as formulating different  related algorithms:
the processes $\PROC{S^m_t}$ corresponding to the  generators $\LOPER_m^E$
are mutually independent  for different $m\in\INDXE$, and thus can be simulated in parallel;  
similarly  we can handle the processes belonging to the
group indexed by $\INDXO$. However, there is still communication between these two groups
as there is non-empty overlap between the groups due to interactions and updates in the sets $\BNDRY\CUBE_m$, $\BNDRY\CUBE_m'$
when $m\in\INDXE$ and $m'\in\INDXO$ and the cells are within the interaction range $L$. 
To handle this communication we next introduce a Fractional Step approximation of 
the Markov semigroup $\EXP{t\LOPER}$ associated with the process $\PROC{S_t}$.  

\subsection{Fractional Step Kinetic Monte Carlo Algorithms}
 The deterministic Trotter Theorem was first proved in \cite{Trotter} for the approximation 
of semigroups corresponding to operator sums, and it has found wide application in the numerical ODE/PDE analysis, e.g., 
\cite{Hairer}. Similarly, the key tool for our analysis is a deterministic as well as 
a stochastic version of the Trotter formula, \cite{Kurtz},
applied to the operator $\LOPER = \LOPER^E + \LOPER^O$
\begin{equation}\label{trotter}
 \EXP{t\LOPER} = \lim_{n\to\infty}\left[ \EXP{\frac{t}{n}\LOPER^E}\EXP{\frac{t}{n}\LOPER^O}\right]^n\PERIOD
\end{equation}
The proposed parallel scheme uses the fact that the action of the operator $\LOPER^E$ (and similarly
of $\LOPER^O$) can be distributed onto independent processing units. Thus to reach a time $T$ we define
a time step $\Delta t = \tfrac{T}{n}$ for a fixed value of $n$ and alternate the evolution by
$\LOPER^E$ and $\LOPER^O$. More precisely,    \VIZ{trotter} gives rise to 
the {\em Lie } splitting approximation for $n\gg 1$:
\begin{equation}\label{lie}
 \EXP{T\LOPER} \approx  \left[ \EXP{\frac{T}{n}\LOPER^E}\EXP{\frac{T}{n}\LOPER^O}\right]^n\PERIOD
\end{equation}  
Since the simulated systems exhibit  short-range interactions, 
the generators  $\LOPER^E_k, \LOPER^E_l$   {\em commute} for $k, l \in\INDXE$, $k\ne l$: 
$$
\LOPER^E_k\LOPER^E_l-\LOPER^E_l \LOPER^E_k=0\, , \quad \mbox{for all $k, l \in\INDXE$, $k\ne l$}\PERIOD
$$
Hence, \cite{Trotter},  we have the exact formula
\begin{equation}\label{exact}
 \EXP{\Delta t\LOPER^E}\EXP{\Delta t\LOPER^O}= \prod_{m\in\INDXE}\EXP{\Delta t \LOPER^E_m}\prod_{m\in\INDXO}\EXP{\Delta t\LOPER^O_m}\PERIOD
\end{equation}
Then the expression \VIZ{exact} implies that  the KMC solvers corresponding to the semigroup 
$\EXP{\Delta t\LOPER^E}$ (resp. $\EXP{\Delta t\LOPER^O}$) can be  simulated  {\em exactly} by breaking down the task
into  separate processors/threads for each $m\in\INDXE$ (resp. $m\in\INDXO$).
Therefore, this  scheme  allows  us to run independently  on each fractional time-step window $\Delta t$,
and on every processor,  a serial KMC simulation, called a {\em kernel\,}. The resulting computational framework consisting of the hierarchical decomposition
\VIZ{opdecomp}  and \VIZ{lie} permits to  input  as  the algorithm's
 kernel any  preferred optimized serial KMC algorithm.

A single time step of the parallel algorithm is thus easily described in the following stages:
\begin{description}
 \item[{\em Step 1--Evolution by $\LOPER^E$}:] Simulate independent Markov processes $\PROC{S_t^m}$, $m\in\INDXE$
                      by a kinetic Monte Carlo kernel running on non-communicating processors that correspond  to 
                       each  $\CUBE_m$  for time  $\Delta t$.
 \item[{\em Step 2--Local Synchronization}:] communicate configurations $\sigma^E$ from overlapping
                     domains $\CLSR\CUBE^E_m\cap\CLSR\CUBE_n^O$ in order to update configurations $\sigma^O$.
 \item[{\em Step 3--Evolution by $\LOPER^O$}:] Simulate independent Markov processes $\PROC{S_t^m}$, 
                     $m\in\INDXO$ by a KMC  {\em kernel}  on non-communicating processors that correspond  to 
                       each  $\CUBE_m$  for time  $\Delta t$.
\end{description}
We emphasize that the resulting process $\PROC{\tilde S_t}$ is 
an {\it approximation} of the process $\PROC{S_t}$ and we discuss it's features and properties  in the   next two sections.

\section{Processor Communication Schedule and Random Trotter Products}\label{pcs}

A key feature of the fractional step methods is the {\em Processor Communication Schedule} (PCS) that dictates the order with which 
the hierarchy of operators in \VIZ{gendecomp} are applied and for how long.  For instance, in \VIZ{lie} the  processors corresponding to $\LOPER^E$ (resp. $\LOPER^O$) do not communicate, hence the processor communication within the algorithm 
occurs {\em only} each time we have to apply   $\EXP{\frac{T}{2n}\LOPER^E}$ or $\EXP{\frac{T}{2n}\LOPER^O}$. 
Therefore we can define as the PCS the (deterministic) jump process $X=X(t)$, $t \in [0, T]$, where $[0, T]$ is the simulated time window and taking values
in the set $\mathcal{X}=\{1, 2\}$, where we assign the value $1$ (resp. $2$) to $O$ (resp. $E$):
\begin{align}\label{PCSl}
 X(t)&=1\COMMA  \quad  \frac{2kT}{n}    \le t < \frac{(2k+1)T}{n}\COMMA \\
 X(t)&=2\COMMA \quad   \frac{(2k+1)T}{n}\le t < \frac{(2k+2)T}{n}\PERIOD
 \end{align}
for all $k=0,\dots, n-1$. Processor communication occurs at jump times, while in the remaining time the processors operate independently
and do not communicate.
In an analogous way we can define the PCS  for the {\em Strang} splitting scheme \VIZ{strang},
\begin{equation}\label{strang}
 \EXP{T\LOPER} \approx \left[ \EXP{\frac{T}{2n}\LOPER^E}\EXP{\frac{T}{n}\LOPER^O}\EXP{\frac{T}{2n}\LOPER^E}\right]^n\COMMA
\end{equation}           
with the scheduling process
\begin{align}\label{PCSs}
  X(t)&=1\COMMA \quad  \frac{2kT}{2n}    \le t < \frac{(2k+1)T}{2n}\COMMA \\
  X(t)&=2\COMMA \quad  \frac{(2k+1)T}{2n}\le t < \frac{(2k+3)T}{2n}\COMMA\\
  X(t)&=1\COMMA \quad  \frac{(2k+3)T}{2n}\le t < \frac{(2k+4)T}{2n}\COMMA
\end{align}
for all $k=0,\dots, n-1$.

\subsection{Random Fractional Step Methods} In both cases above \VIZ{lie} and \VIZ{strang}, the communication schedule 
is fully deterministic, relying on the Trotter Theorem \VIZ{trotter}.
On the other hand, we can construct stochastic PCS based on the  {\em Random Trotter Product} Theorem, and 
as we show below the 
sub-lattice algorithm proposed in \cite{ShimAmar05b} is a fractional step algorithm with stochastic PCS. 
    
The   Random Trotter Product Theorem, \cite{Kurtz},  extends \VIZ{trotter} as follows:
We consider a sequence of semigroups $ \EXP{T\LOPER_\xi}$ with corresponding operators $\LOPER_\xi$ where $\xi$ is  in the index set  
$\mathcal{X}$, assuming for simplicity $\mathcal{X}$  is finite, although a much more general setting is possible, \VIZ{kurtz}.  
Consider also a stochastic jump process $X=X(t)$ with $\mathcal{X}$ as its state space. 
For each of its trajectories we denote by $\xi_0, \xi_1, ...\xi_n$ the (typically random) sequence of states visited by the stochastic process $X(t)$ 
and $\tau_0, \tau_1, \dots, \tau_n$ the corresponding (also typically random) jump times
\begin{align}\label{PCSr}
  X(t)&=\xi_0\COMMA \quad  0        \le t< \tau_0\COMMA\\
  X(t)&=\xi_1\COMMA \quad  \tau_0   \le t< \tau_1\COMMA\\
      ... \\
  X(t)&=\xi_k\COMMA \quad \tau_{k-1}\le t< \tau_k\PERIOD
\end{align} 
 We additionally define as $N(t)$ the number of jumps up to time $t$. We assume that  $X(t)$ is selected so 
 that it has an ergodic behavior, i.e., there is a probability measure $\mu(d\xi)$ such that  for all bounded  
functions $g$ we have that
 \begin{equation}\label{ergodic}
   \lim_{t \to \infty} \frac{1}{t} \int_0^t g(X(s))\, ds = \int g(\xi)\, \mu(d\xi)\PERIOD
 \end{equation}
 For example, if $X(t)$ is a Markov process then under suitable conditions, \VIZ{ergodic} will hold, where $\mu$ will be the stationary distribution
 of $X(t)$, \cite{Liggett}.
 Conversely, it is well-known that for a given $\mu$ we can construct in a non-unique way Markov processes  $X(t)$ which satisfy  
the condition \VIZ{ergodic}, \cite{Liggett}.
 Now we can state the Random Trotter Product Theorem, \cite{Kurtz}, in analogy to \VIZ{trotter}:
\begin{equation}\label{kurtz}
   \EXP{T\bar\LOPER} = \lim_{n\to\infty}\left[ \EXP{{ \tau_0\over n} \LOPER_{\xi_0}}\EXP{{\tau_1-\tau_0 \over n}\LOPER_{\xi_1}}\dots\EXP{{nT-\tau_{N(nt)}\over n}\LOPER_{\xi_{N(nt)}}}\right]\COMMA
\end{equation} 
where the operator $\bar\LOPER$ is defined on any bounded function as
\begin{equation}\label{averaged}
   \bar\LOPER g=\int \LOPER_\xi \mu(d\xi)\, .
\end{equation}
It is clear that \VIZ{lie} is a special case of \VIZ{kurtz} when $\tau_k-\tau_{k-1}=1$ and 
$\xi_{2k}=1$, $\xi_{2k+1}=2$ for all $k$. Similarly, we can also view  \VIZ{strang} as a deterministic analogue of \VIZ{kurtz}. 

On the other hand, in the context of the parallel fractional step algorithms for KMC introduced here,
the random process \VIZ{PCSr} can be interpreted as a stochastic PCS.  For example,  the  sub-lattice (SL) 
parallelization algorithm for KMC,
introduced in \cite{ShimAmar05b},  is a fractional step algorithm with stochastic PCS: indeed, 
in this method the lattice is divided into  sub-lattices, for instance as in \VIZ{sublatt},
$\LATT = \LATT^E\cup\LATT^O$. Each sub-lattice is selected {\em at random} and  advanced by KMC over a fixed time window $\Delta t$.
Then a new random selection is made and again the sub-lattice is advanced by $\Delta t$, and so on. The procedure is parallelizable as cells $\CUBE_m^E$,
$\CUBE_m^O$ within each sub-lattice do not communicate. 
This algorithm is easily recast as a fractional step approximation, when in \VIZ{PCSr} we select deterministic jump times  
$\tau_k$ and random variables $\xi_k$:
\begin{equation}\label{SLPCS}
  {\tau_k-\tau_{k-1}\over n}=\Delta t\COMMA \quad \mbox{and}\quad  P(\xi_k=1)=P(\xi_k=2)=\frac{1}{2}\PERIOD
\end{equation}
As in \VIZ{PCSl}, here we assign the value $1$ (resp. $2$) to the $O$ (resp. $E$) sub-lattice. Furthermore, we can easily calculate \VIZ{averaged}
to obtain 
$$
 \bar\LOPER g=\frac{1}{2} \left(\LOPER^E+\LOPER^O\right)\COMMA
$$
which is just a time rescaling of the original operator $\LOPER$. Thus the SL algorithm is rewritten as  the fractional step approximation 
with the  stochastic PCS \VIZ{SLPCS} as 
\begin{equation}\label{SL}
 \EXP{T\bar\LOPER} \approx  \EXP{{ \tau_0\over n} \LOPER_{\xi_0}}\EXP{{\tau_1-\tau_0 \over n}\LOPER_{\xi_1}}\dots\EXP{{nT-\tau_{N(nt)}\over n}\LOPER_{\xi_{N(nt)}}} \PERIOD
\end{equation}
From the numerical analysis viewpoint, our   re-interpretation  of the SL algorithm in  \cite{ShimAmar05b} as a fractional step scheme allows us 
to  also provide a mathematically  rigorous justification  that it is a {\em consistent} estimator of the serial KMC algorithm,  due to  the Random Trotter Theorem \VIZ{kurtz}. 
That is,  as the time step in the fractional step scheme converges to zero,  
it  converges to the continuous time Markov Chain that has the same master equation and generator as 
the original serial KMC. Finally,  the (deterministic) Trotter Theorem \VIZ{trotter} also implies that the 
Lie and the Strang schemes are, in the numerical analysis sense,  consistent  approximations of the serial KMC algorithm.%

\section{Controlled Error Approximations of KMC}\label{cerror}

In this section we present  a formal argument for  the error analysis of the fractional step approximations  for KMC, which suggests the order of convergence of the schemes, as well as the restrictions on the Fractional Step KMC time step $\Delta t$.%
In the decomposition \VIZ{opdecomp} the operators  
are linear operators on the high, but finite-dimensional configuration space $\SIGMA$, hence  by the standard error analysis 
of splitting schemes, see \cite{Hairer}, we have
\begin{equation}\label{error}
\EXP{\Delta t \LOPER} -\EXP{\Delta t \LOPER^E}\EXP{\Delta t \LOPER^O}=[\LOPER^E\LOPER^O-\LOPER^O\LOPER^E]{(\Delta t)^2\over 2}+\BIGO(\Delta t^3)\COMMA
\end{equation}
where we  readily see that the term  $[\LOPER^E, \LOPER^O]:=\LOPER^E\LOPER^O-\LOPER^O\LOPER^E$ is the Lie bracket (commutator) of the operators 
$\LOPER^E$, $\LOPER^O$. This Lie bracket captures the effect of the
boundary regions $\CLSR\CUBE^E_m\cap\CLSR\CUBE_n^O$ through which we have processor communication: if there was no communication the Lie bracket would be exactly zero.

Furthermore, instead of \VIZ{lie} we can consider  the {\em Strang-type} splitting \VIZ{strang}.
As in the ODE  case,  \cite{Hairer}, this  is expected to yield a higher order error term $\BIGO(\Delta t^3)$ instead of the second order approximation in \VIZ{error}, in the following sense:
\begin{align}\label{error2}
\EXP{\Delta t \LOPER} -\EXP{\frac{\Delta t}{2}\LOPER^E}\EXP{\Delta t\LOPER^O}\EXP{\frac{\Delta t}{2}\LOPER^E}
= \, &\Big\{{1 \over 12}[\LOPER^O, [\LOPER^O, \LOPER^E]\, ]\nonumber\\
-{1 \over 24}[\LOPER^E, [\LOPER^E, &\LOPER^O]\, ] \Big\}(\Delta t)^3+\BIGO(\Delta t^4)\PERIOD
\end{align}
Such calculations  suggest that   the Strang splitting  leads to a more accurate scheme, which is  balanced by   
more complicated boundary local communication  in the same time window $\Delta t$, as is evident 
when comparing \VIZ{lie} and  \VIZ{strang}.

Next,  we briefly comment on the error estimation suggested by the calculation  \VIZ{error} and return to the rigorous numerical analysis in \cite{AKPR11}. 
In order  to obtain an  estimate in the right-hand side of \VIZ{error} which  is  independent of the system size $N$,  it is essential to obtain an upper bound on the total number 
of jumps up to the time $T$. This is a key point related to the {\em extensivity} of the system and to the fact that the  weak error analysis is  restricted (as it should be physically)
to mesoscopic observables satisfying \VIZ{observables}. We observe the dependence of the error on mesoscopic observables in the following subsection.
In the context of coarse-graining, in \cite{KPS} an analogous estimate   was shown rigorously using a Bernstein-type argument 
applied to the 
discrete derivatives, in the spirit of \VIZ{observables},  of the solutions to the backward Kolmogorov equation. We refer to  such bounds as ``Bernstein-like'' due to their  similarity  to gradient  estimates for linear 
and nonlinear parabolic PDEs.

\medskip

\subsection{Error Analysis and comparison between random and deterministic PCS}\label{comparison}
In this section we further demonstrate the use of the operator splitting formulation as a numerical analysis tool by comparing the time-step of $\Delta t$  the random PCS introduced 
in \cite{ShimAmar05b} to the deterministic Lie PCS introduced in \VIZ{lie}. A similar comparison can be made for the Strang scheme \VIZ{strang}.
A detailed discussion including rigorous error estimates for mesoscopic observables such as \VIZ{observables}, which are independent of the lattice size
$N$ will be discussed in \cite{AKPR11}.

Here we focus on the example of adsorption/desorption discussed in Section \ref{formulation}. The generator in the one space dimension is decomposed
as in \eqref{opdecomp}
\[
\mathcal{L}^Ef(\sigma) = \sum_{x\in\Lambda} c^E(x,\sigma) \Bigl( f(\sigma^x)-f(\sigma) \Bigr)\COMMA
\]
and
\[
\mathcal{L}^Of(\sigma) = \sum_{x\in\Lambda}  c^O(x,\sigma) \Bigl( f(\sigma^x)-f(\sigma) \Bigr)\COMMA
\]
where
\[
  c^E(x,\sigma) = \left\{ \begin{array}{ll} c(x,\sigma), & x\in \Lambda^E_N \\ 0 , &  \mbox{otherwise} \end{array} \right . \;\;\;\;
  c^O(x,\sigma) = \left\{ \begin{array}{ll} c(x,\sigma), & x\in \Lambda^O_N \\ 0 , & \mbox{otherwise} \end{array} \right .
\]
and the sub-lattices $\Lambda_N^E, \Lambda_N^O$ are defined in \eqref{sublatt}. The rates $c(x, \sigma)$ of the corresponding generator \VIZ{generator}
for the case of Arrhenius adsorption/desorption are given by
\begin{equation}\label{adsdesrate}
c(x, \sigma)=c_a(1-\sigma(x))+c_d\sigma(x)\exp{\big(-\beta U(x, \sigma)\big)}\, ,
\end{equation}
where $c_a$ and $c_d$ are the adsorption and desorption constants respectively, \cite{Vlachos}. The desorption potential $U=U(x, \sigma)$ is defined as
\begin{equation}\label{potential}
U(x, \sigma)=\sum_{y\ne x}J(x-y)\sigma(y)\, ,
\end{equation}
where $J=J(x-y)$ is the lateral interaction potential; for simplicity we assume that the range of interactions is $L$, while in typical
simplified nearest neighbor models $L=1$. Similarly we define diffusion dynamics with Arrhenius dynamics, \cite{KMV}.

First we discuss the error analysis for the Lie splitting scheme. For given  finite lattice size $N$,
in the decomposition \VIZ{opdecomp} the operators  
are linear operators on the high, but finite-dimensional configuration space $\SIGMA$, hence  by the standard error analysis 
of Lie splitting schemes, we obtain \VIZ{error}.
A more careful study of the commutator reveals that the generator decomposition \VIZ{opdecomp}
induces significant cancellations in the evaluation of the generator: indeed, we define
\[
\CUBE_m^o=\CUBE_m \setminus  \BNDRY\CUBE_m\, , \quad\quad \CUBE_m=\CUBE_m^o \cup \CUBE_m^{\partial}\, ,
\]
where  in Section~\ref{formulation} we introduced
$\BNDRY\CUBE_m = \CLSR\CUBE_m\cap\CUBE_m$ and 
$\CLSR\CUBE_m = \{ z\in\LATT\SEP |z - x|\leq L\COMMA x\in\CUBE_m\}\PERIOD
$
Thus, in \VIZ{opdecomp} we obtain the further decomposition 
\begin{equation}\label{opdecomperror}
\LOPER^E=\LOPER^{E, o}+\LOPER^{E, \partial}:=\sum_{m \in \INDXE}\LOPER_m^{E, o}+\LOPER_m^{E, \partial}\, ,
\end{equation}
where $\LOPER_m^{E, o}, \LOPER_m^{E, \partial}$ is the restriction of $\LOPER^E$ on $\CUBE_m^o$ and $\CUBE_m^{\partial}$ respectively.
Analogously we define $\LOPER^O=\LOPER^{O, o}+\LOPER^{O, \partial}$.
We now return to the evaluation of the commutator
\begin{equation}
 [\mathcal{L}^E, \mathcal{L}^O ] = [ \mathcal{L}^{E,\partial} , \mathcal{L}^{O,\partial}]+[ \mathcal{L}^{E,o} , \mathcal{L}^{O,o}]+ 
 [\mathcal{L}^{E,\partial},\mathcal{L}^{O,o}] +[ \mathcal{L}^{E,o} , \mathcal{L}^{O,\partial}] \PERIOD
 \end{equation}
However, due to the lack of communication between generators beyond the interaction range, we have that 
\[
 [ \mathcal{L}^{E,o} , \mathcal{L}^{O,o}]=0\, , \quad
 [ \mathcal{L}^{E,\partial} , \mathcal{L}^{O,o}]=0\, , \quad [ \mathcal{L}^{E,o} , \mathcal{L}^{O,\partial}]=0\, ,
\] 
thus we readily get
\begin{equation}
\label{lieboundary}
 [\mathcal{L}^E, \mathcal{L}^O ] = [ \mathcal{L}^{E,\partial} , \mathcal{L}^{O,\partial}]=
 \sum_{m \in \INDXE}
 \sum_{\substack{l\in\INDXO \\ 
 |l-m|=1}}
[ \mathcal{L}_m^{E,\partial} , \mathcal{L}_{l}^{O,\partial} ]
\, . 
 \end{equation}
The formula \VIZ{lieboundary}  captures the processor communication between 
boundary regions of $\CLSR\CUBE^E_m$, $\CLSR\CUBE_n^O$. But more importantly, when combined with  \VIZ{error}, it suggests  
the limitations on the time window $\Delta t$ of the Lie scheme \VIZ{lie}, denoted for differentiation by $\Delta t_{\mbox{Lie}}$,  in order to obtain a given error tolerance TOL.
In that sense it is useful to obtain an upper bound on \VIZ{lieboundary}. Indeed, we readily obtain:
\begin{align}
[\mathcal{L}^E , \mathcal{L}^O ]f(\sigma) = 
 \sum_{\substack{m \in \INDXE, l\in\INDXO \\ 
 |l-m|=1}}&\sum_{x,y}\Bigl [ c^E(x,\sigma)c^O(y,\sigma^x) - c^E(x,\sigma^y)c^O(y,\sigma) \Bigr ] f\big((\sigma ^x)^{y}\big) \nonumber\\
& -\sum_{x, y} c^E(x,\sigma) \Bigl [ c^O(y,\sigma^x) - c^O(y,\sigma) \Bigr ] f(\sigma ^{x})  \nonumber\\
 &-\sum_{x,y} c^O(y,\sigma) \Bigl [  c^E(x,\sigma) -  c^E(x,\sigma^y) \Bigr ] f(\sigma ^{y})
\end{align}
where all  summations are over $x \in C_m^{E, \partial},y \in C_l^{O, \partial}$.
For    {\em mesoscopic observables}, such as the mean coverage  $f(\sigma) = {1 \over N} \sum_{x\in\Lambda}\sigma(x)$ we obtain
\begin{align}
[ \mathcal{L}^E , \mathcal{L}^O ] f(\sigma) =  \sum_{\substack{m \in \INDXE, l\in\INDXO \\ 
 |l-m|=1}}&\sum_{x,y} c^O(y,\sigma) \Bigl [  c^E(x,\sigma) -  c^E(x,\sigma^y) \Bigr ] { 1- 2 \sigma(x)\over N}  \nonumber\\
&+\sum_{x,y}  c^E(x,\sigma) \Bigl [ c^O(y,\sigma^x) - c^O(y,\sigma) \Bigr ] {1-2\sigma(y)\over N}\, ,
\end{align}
where all  summations are over $x \in C_m^{E, \partial},y \in C_l^{O, \partial}$.
Therefore, due to the {\em cancellation}
of all interior components $\LOPER^{E, o}, \LOPER^{O, o}$ in \VIZ{lieboundary}, we obtain the bound for the case of  the interaction range $L=1$,
\begin{equation}
\label{liebound}
|[ \mathcal{L}^E , \mathcal{L}^O ] f(\sigma) | \sim O\Big({M\cdot L \over N}\Big)=O\Big({1 \over q}\Big)\, ,
\end{equation}
where $q$ is the size of each cell $\CUBE_m$, and $O(1)$ depends on the physical parameters in the rate \VIZ{adsdesrate}.
The local error analysis in \VIZ{error}, \VIZ{liebound}  can be propagated up to a prescribed time $T=N_{\mbox{Lie}}\Delta t_{\mbox{Lie}}$
Therefore, for the simulation of the mesoscopic observable $f$ up to the time $T$ within a given error tolerance TOL, \VIZ{error} and \VIZ{liebound} give the 
{\em observable}-dependent relation for the Lie time step
\begin{equation}\label{lietimestep}
\mbox{TOL}\sim T \cdot  |[ \mathcal{L}^E , \mathcal{L}^O ] f(\sigma) |\Delta t_{\mbox{Lie}} \sim T\cdot O\Big({1 \over q}\Big)\Delta t_{\mbox{Lie}}
\end{equation}

Next, using the fractional step formulation, we analyze in the same spirit as for the Lie scheme, the random PCS \VIZ{SLPCS} proposed in \cite{ShimAmar05b}.
For notational simplicity we set $A_1=\mathcal{L}^O$, $A_2=\mathcal{L}^E$. Then the local error operator $E^{\Delta t}$ can also be calculated as in \VIZ{error}:
\begin{align}
\mbox{Local Error}=E^{\Delta t}:=& e^{\Delta t A_{\xi_1}}e^{\Delta t A_{\xi_2}} - e^{\Delta t (A_1+A_2)} \nonumber\\
=& \Bigl ( I + (A_{\xi_1}+A_{\xi_2})\Delta t + \frac{1}{2}(A_{\xi_1}^2+2A_{\xi_1}A_{\xi_2} +A_{\xi_2}^2)\Delta t^2 \Bigr ) -\nonumber \\
&\Bigl ( I + (A_{1}+A_{2})\Delta t + \frac{1}{2}(A_{1}+A_{2})^2\Delta t^2 \Bigr ) + O(\Delta t^3)
\end{align}
The mean value of the error over the sequence of independent random variables $\xi=(\xi_i\, , i=1, ..., n)$ of the PCS \VIZ{SLPCS} on an observable $f=f(\sigma)\, , s \in \SIGMA$ can be  explicitly evaluated:
\[
\EXPECT_{\xi}[E^{\Delta t}f] = \frac{1}{4}(A_1-A_2)^2f\Delta t^2+ O(\Delta t^3)\quad= \frac{1}{4}(\LOPER^E-\LOPER^O)^2f\Delta t^2+ O(\Delta t^3)\, .
\]
As in \VIZ{liebound}, for  the   mesoscopic observable  $f(\sigma) = {1 \over N} \sum_{x\in\Lambda}\sigma(x)$, we obtain, after disregarding the higher order local error $O(\Delta t^3)$, 
\begin{equation}\label{randombound}
(\LOPER^E-\LOPER^O)^2f(\sigma)\sim O(1)\, ,
\end{equation}
where $O(1)$ depends on the physical parameters in the rate \VIZ{adsdesrate}. Similarly to \VIZ{lietimestep}, for the simulation of the mesoscopic observable $f$ up to the same prescribed time $T=N_{\mbox{Random}}\Delta t_{\mbox{Random}}$, within the same error tolerance TOL, \VIZ{error} and \VIZ{randombound} give the 
{\em observable}-dependent relation for the random PCS  time step
\begin{equation}\label{randomtimestep}
\mbox{TOL}\sim  T\cdot |(\LOPER^E-\LOPER^O)^2f(\sigma) |\Delta t_{\mbox{Random}} \sim T \cdot O(1)\Delta t_{\mbox{Random}}
\end{equation}
Comparing the random and the Lie PCS through \VIZ{lietimestep} and \VIZ{randomtimestep} implies that in order the two schemes to conform (in the mean) to the same tolerance TOL, their respective time steps should be selected so that 
\begin{equation}\label{pcscompare}
 \Delta t_{\mbox{Lie}} \sim O(q)\Delta t_{\mbox{Random}}
\end{equation}
This relation in turn suggests that  the Lie scheme \VIZ{lie} is expected to parallelize better than the random PCS \VIZ{SLPCS} since it allows a $q$-times larger time step $\Delta t$ for the same accuracy, keeping in mind that  during each time step  processors do not communicate.

A similar analysis is possible for general mesoscopic  observables $f=f(\sigma)\, , s \in \SIGMA$, e.g., spatial correlations, that satisfy
\begin{equation}\label{observables}
\sum_{x \in \LATT}|f(\sigma^x)-f(\sigma)| \le C 
\end{equation}
where $C$ is a constant independent of $N$, see the formulation and estimates for coarse-grained stochastic systems in \cite{KPS}. We revisit this issue, as well as the rigorous derivation of $N$-independent error bounds
in place of  the expansions \VIZ{error}, \VIZ{error2} in the upcoming publication \cite{AKPR11}. Such estimates can also allow a  detailed analysis on the balance between accuracy and local processor communication for PCS such as \VIZ{lie}, \VIZ{strang} and \VIZ{SLPCS}.

 \section{Hierarchical structure of Fractional Step algorithms and implementation on  GPUs}\label{gpu}
 The fractional step framework  allows a hierarchical structure to be easily formulated and implemented, which is a key  advantage for simulating 
in  parallel architectures with complex  memory hierarchies and processing units.
The Graphical Processing Unit (GPU) architecture is inherently
different from a traditional CPU architecture. GPUs are massively
parallel multi-threaded devices capable of executing a large number of
active threads concurrently. A GPU consists of multiple streaming
multiprocessors (MP), each of which contains multiple scalar processor
cores. For example, NVIDIA's C2050 GPU architecture contains 14 such
multiprocessors, each of which contains 32 cores, for a total of 448
cores which can handle up to 24k active threads in parallel. A GPU has
several types of memory which are differently organized compared to
the traditional hierarchical CPU memory, most notably the main device
memory (global memory) shared between all the multiprocessors and the
on-chip memory shared between all cores of a single multiprocessor
(shared memory). The memory sizes and access speeds depend on the type
of GPU. For instance, the memory size of the NVIDIA C2050 GPU is 3GB
while the memory size of the NVIDIA C2070 GPU is 6GB.

From the perspective of a GPU programmer writing a code for NVIDIA
GPU's, the GPU is treated as a co-processor to the main CPU. Programs
are written in C and linked to the CUDA libraries. A function that
executes on the GPU, called a GPU kernel, consists of multiple threads
executing code in a single instruction, multiple data (SIMD)
fashion. That is, each thread in a GPU kernel executes the same code,
but on different data. 
Further, threads can be grouped into
thread blocks. This abstraction takes advantage of the fact that
threads executing on the same multiprocessor can share data via
on-chip shared memory, allowing some degree of cooperation between
threads in the same block~\cite{CUDA}. A major drawback in GPU
programming is the slow communication between GPU global memory and
the main memory of the CPU, compared to the communication within a
GPU. Programmers address this problem by maximizing the amount of
arithmetic intensive computations performed on GPU, minimizing the
communication between CPU and GPU, and allowing the communication
latency to be hidden by overlapping with execution. Communication
among GPUs, although costly, is enabled by APIs such as OpenMP and
features available in CUDA 2.2+ such as portable pinned memory, when
the communication is among GPUs connected to the same shared-memory
computer node. When the communication takes place among GPUs across
nodes of a cluster, message passing paradigms such as MPI can serve
the same scope.
\begin{figure}
\centerline{\includegraphics[width=.8\textwidth]{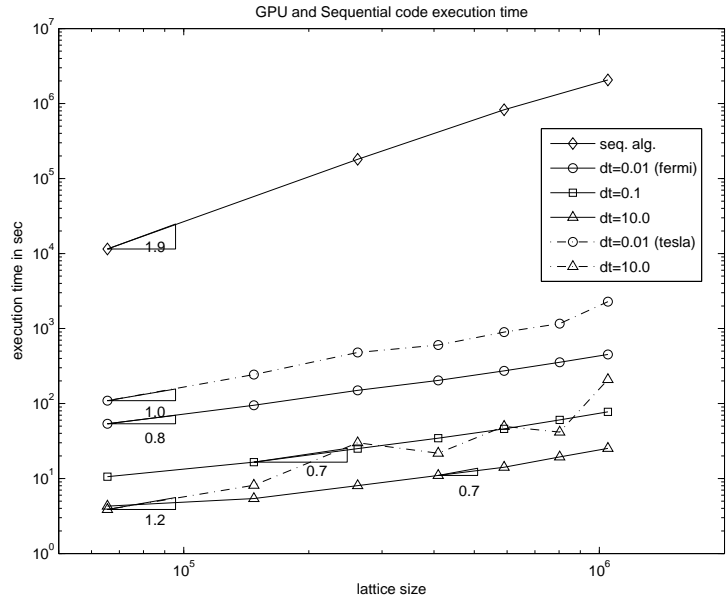}}
\caption{Execution time of the fractional step KMC  for lattices of different sizes. The comparison with the sequential algorithm (top curve)
         is based on the same SSA KMC implementation which, however, does not have the optimal complexity of the BKL algorithm.
         The  simpler implementation
         of the SSA algorithm was used. The simple implementation has the complexity $\BIGO(N^2)$, where $N$ is the total number of lattice
         sites. This complexity is reflected in the indicated scaling (the slope in the log-log plot).
         Note that the due to partitioning of the lattice in the fractional step algorithm the same KMC kernel
         will scale as $\BIGO(N)$ only, which is in agreement with the observed slope in the plots.}\label{timeGPUs}
\end{figure}

In our parallelization of the KMC method, we redefine the data
structures to represent lattice sites in the simulation so that the
whole simulated system is cut into equal-sized black and white coarse cells
like a chessboard in \VIZ{sublatt}. For instance, 
Fig.~\ref{latt1}
shows a simple example in which we map a $4\times 4$ lattice sites into
$2\times 2$ cells, each cell containing $2\times 2$ sites. One GPU thread is
assigned to one cell. Coverage information of the whole lattice is
stored in an array located in the GPU global memory so that all the
threads can access the information related to their neighboring sites
across MPs. The GPU kernel performing the KMC simulation over the whole lattice by
using the Lie scheme \VIZ{lie} and 
 the  decomposition \VIZ{opdecomp}, is
sequentially launched twice for each synchronization time step $\Delta t$
to work on the black and white cells respectively.
 The execution times for lattices of different sizes are compared in Fig~\ref{timeGPUs},
 where we take as a reference a sequential KMC-kernel, which is a direct numerical implementation of
 \VIZ{totalrate} and \VIZ{skeleton}. The same kernel is then used
for the implementation on GPUs where we compare times for different choices of $\Delta t$. 
We remark that the KMC kernel is not optimized by techniques such as the BKL algorithm, 
\cite{BKL75, binder},
which is also manifested in the scaling with respect to the size of the lattice $N$. However, the same
kernel is used in the Fractional Step algorithm thus  here we present comparisons between the same KMC algorithms,
one serial and one parallelized by the Fractional Step approach. Clearly any optimized  KMC kernel can be used
without difficulty in our framework. 

The size of lattices that can be simulated on a single GPU is limited by  memory,  thus in
order to simulate large systems  it will be necessary  to employ a cluster of GPUs communicating, for instance, through an MPI protocol.
We will demonstrate next how Fractional Step KMC   algorithms can 
be tailored to  an  architecture that involves multiple GPUs.
We  return  to the formulation in \VIZ{opdecomp}, and consider the sub-lattice
decomposition \VIZ{sublatt}. In this formulation each one of the coarse-cells $\CUBE_m^E$ or $\CUBE_m^O$ are simulated on a single GPU. Within each one of the GPUs  we have the same lattice decomposition as in \VIZ{sublatt}, see Figure~\ref{latt2}, namely
\begin{equation}\label{sublatt2}
\CUBE^E_m= \CUBE_m^{EE}\cup\CUBE_m^{EO}:=\bigcup_{l=1}^L D_{ml}^{EE} \cup \bigcup_{l=1}^L D_{ml}^{EO}  \COMMA
\end{equation}
and similarly we define a decomposition  for $\CUBE^O_m$.
Each one of the (sub-)sub-lattices $D_{ml}^{EE}$ and $D_{ml}^{EE}$ corresponds to individual  threads within the GPU. Next,  \VIZ{sublatt} and \VIZ{sublatt2} define  {\em nested sub-lattices}, 
which yield a hierarchical decomposition of the operator $\LOPER$ into \VIZ{opdecomp} and 
\begin{equation}\label{opdecomp2}
\LOPER^E_m=\LOPER_m^{EE}+\LOPER_m^{EO}:= \sum_{l=1}^L \LOPER_{ml}^{EE} + 
\sum_{l=1}^L \LOPER_{ml}^{EO} 
\COMMA
\end{equation}
and similarly we also define the decomposition for $\LOPER^O_m$. 
Finally, schemes such as  \VIZ{lie} and \VIZ{strang} give rise to Fractional Step algorithms based on the nested decompositions \VIZ{opdecomp} and \VIZ{opdecomp2}. 
In this case, boundary communication, see Fig.~\ref{latt2}, plays a key role in the parallelization of our
algorithm when multiple GPUs are required. As we discussed earlier, this scenario happens when
the lattice size grows to the point that the lattice data structures
no longer fit into a single GPU global memory. In turn, this threshold
depends on the type of GPU used, e.g., for a NVIDIA's C2050 GPU the
maximum lattice size is currently $8,182\times 8,182$ cells. To simulate larger
systems, we can decompose the domain into regular sub-domains and
distribute both the sub-domain cells and associated computation among
multiple GPUs, as discussed in \VIZ{opdecomp2}. Boundary communication between two adjacent sub-domains
are exchanged between GPUs, see Fig.~\ref{latt2}, and supported by either MPI or OpenMP,
depending on the fact that the GPUs are located on the same cluster
node or across nodes. Thus, the multi-GPU parallel KMC algorithm is based 
on and benefits from the hierarchical structure of the
Fractional Step KMC algorithms discussed in \VIZ{opdecomp2}. At
the same time, it can enable the scalability of our simulations to lattice sizes beyond
the ones accessible with a single GPU
e.g.,  $8,182 \times 8,182$ sites in a C2050 GPU.
The study of performance and
scalability of our multi-GPU algorithm and code for different lattice sizes and
types of GPU clusters is beyond the scope of this paper.

\section{Mass Transport and Dynamic Workload Balancing}\label{workloadb} 
 
\begin{figure}
\centerline{
    \subfigure[]{\label{snapshot1}\includegraphics[width=0.55\textwidth]{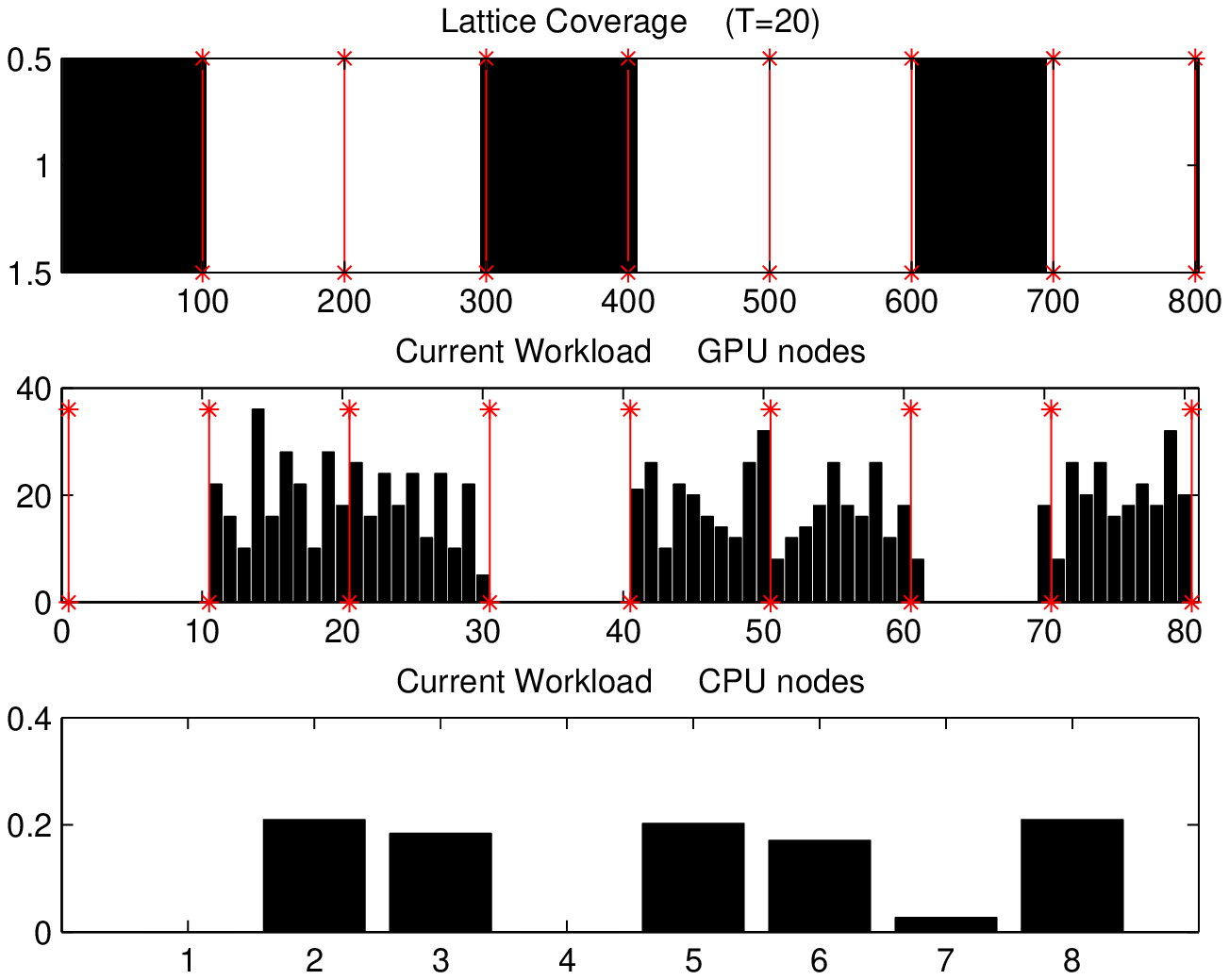}}
    \subfigure[]{\label{snapshot2}\includegraphics[width=0.55\textwidth]{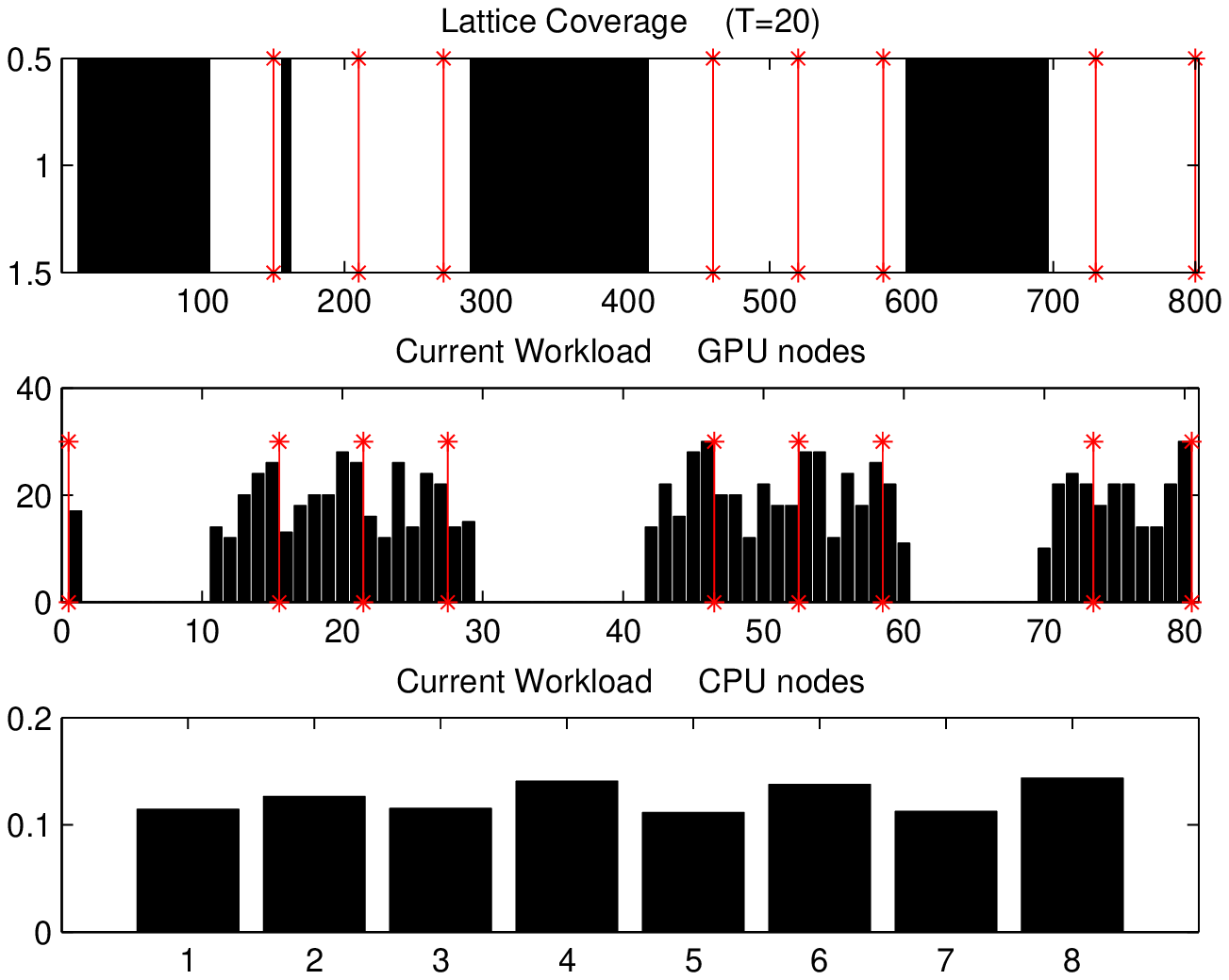}}
  }
\caption{(a) Workload imbalance in 1D unimolecular reaction system: the top figure depicts local coverage, the bottom  figure workload distribution; 
         (b) Workload redistribution in Figure (a)  using the mass transport for re-balancing. }\label{balance}
\end{figure}

Due to the spatially distributed nature of  KMC simulations and the dependence of jump rates on local coverage, \VIZ{totalrate},
fractional step algorithms may have an imbalance in the number of operations/jumps  performed in each coarse cell $\CUBE_m$ in \VIZ{sublatt}, 
as well as on  the corresponding processors. In fact,  formulas \VIZ{totalrate} and \VIZ{skeleton}, and the very structure  of the 
fractional step algorithms \VIZ{opdecomp}, allow us to define the {\em workload} $W_{n\Delta t}(\sigma)=W_{n\Delta t}(m; \sigma)$, $1\le m\le M$ as
\begin{equation}\label{wload}
  W_{n\Delta t}(m)=\mbox{$\#$ jumps in $\CUBE_m$ during  $\left[(n-1)\Delta t, n \Delta t\right]$}\COMMA
\end{equation}
when the configuration at time $(n-1)\Delta t$ is $\sigma$. We also  renormalize $W_{n\Delta t}$ (and still denote it with the same symbol) 
in order to obtain a histogram, i.e., a probability density. Since different coarse cells $\CUBE_m$ in the fractional step algorithms such as \VIZ{lie}
or \VIZ{strang} do not communicate during intervals of length  $\Delta t$ the quantities $\VIZ{wload}$ are easy to keep track on-the-fly during the simulations.
The  possibility of workload imbalance is depicted in Figure~\ref{balance}, where many more jumps are performed in  the processors
corresponding to cells of low coverage, while the other processors remain idle.
 
In this Section we introduce  a probabilistic strategy to re-balance the workload  $W_{n\Delta t}$ dynamically during the simulation 
based on the following  idea from  {\em Mass Transport} methods, e.g., \cite{lcevans}. 
One wants to transport the ``imbalanced'' density $W_{n\Delta t}$ into an almost uniform density over the number of processors used, 
in order to ensure that they remain as uniformly active as possible. The mass transport connection and terminology refers to the mapping  
of a given probability measure into a desirable probability measure. Typically, \cite{lcevans}, this problem is posed as 
an optimization over a suitable cost functional and is known as the Monge-Kantorovich problem. 
In our context the cost functional could reflect constraints related to various parallel  architectures. 

We can formulate and implement this strategy in several different ways:
probably the simplest approach, that serves mostly as an  illustration,  is to assume that we have a number of processors $P$, where $P \ll M$; 
during the interval $\left[(n-1)\Delta t, n \Delta t\right]$ a number of coarse cells $\CUBE_m$, $1\le m\le M$, 
which are simulated independently in a fractional step algorithm, are allocated to each processor. By the end of the simulation time $n\Delta t$ 
the workload on all processors is described similarly to \VIZ{wload}, by a histogram $R_{n\Delta t}(\sigma)=R_{n\Delta t}(l; \sigma)$, $1\le l \le P$. 
One wants to map \VIZ{wload} onto a histogram 
$R_{n\Delta t}$ which is almost uniform in $1\le l\le P$. One such  function can be constructed by  mapping the  mass corresponding to each value  
of the cumulative distribution function (cdf)  of  \VIZ{wload},  onto an equal mass
on the uniform distribution over the $P$ processors. 
In another implementation of the mass transport method  we can adjust the size of the coarse cells $\CUBE_m$ according to the workload redistribution 
strategy discussed earlier, see  Figure~\ref{balance}.  This is effectively a one-dimensional example of an adsorption/desorption process
where the mass transport procedure is carried out by mapping \VIZ{wload}  into a new histogram $R_{n\Delta t}(\sigma)=R_{n\Delta t}(l; \sigma)$ 
corresponding to a new set of variable size coarse cells $\CUBE_l$,  $1\le l \le M'$. 
The cell size adjustment ensures the uniformity of the new histogram by defining $R_{n\Delta t}$ as a mapping  of the cdf corresponding to \VIZ{wload}.

The mass transport  mappings discussed above are not expected to    
 be carried out at every time step $n \Delta t$ in order to reduce computational and communication cost,
 but instead they should   follow a rationally designed coarser-in-time schedule, in analogy to  processor communication scheduling, e.g., \VIZ{PCSr}.The overall implementation appears rather simple since here we demonstrated the methodology 
in a one-dimensional example. However, in higher dimensions, adjusting the size and shape of coarse cells $\CUBE_m$ can be much harder.
Nevertheless the structure of re-balancing procedure can remain one-dimensional even in higher dimensional lattices  if we pick a sub-lattice decomposition \VIZ{decomposition} into strips $\CUBE_m$.
We note that   the mapping we constructed using cdf's did not take into account the processor architecture and a suitable cost functional 
formulation for the mass transport to a uniform distribution, as in the Monge-Kantorovich problem, \cite{lcevans}, may be more appropriate.  We will revisit such issues in a future publication.

\section{Parallel Simulations: Benchmarks and Applications}\label{apps}
Exactly solvable models of statistical mechanics provide a test bed for sampling algorithms applied to
interacting particle systems. We present benchmarks for two important cases: 
(a) sampling of equilibrium distributions,
i.e., long time behavior of the simulated Markov process, and 
(b) weak approximations of the dynamics.
In the first set of tests we work with the classical Ising model on one and two dimensional lattices
where spins interact through a nearest-neighbor potential. Thus the Hamiltonian of the system is
$$
H(\sigma) = -\frac{K}{2}\sum_{x\in\LATT}\sum_{|y-x|=1} \sigma(x)\sigma(y) + h\sum_{x\in\LATT}\sigma(x)\COMMA
$$
where $K$ is a real parameter that defines the strength of the interaction and $h$ the external field.
We work with the spin-flip Arrhenius dynamics with the rates defined in the
nearest-neighbor set $\Omega_x= \{z\SEP |z-x|=1\}$ and the updates in $\SIGMA_x = \{0,1\}$.
\begin{eqnarray}\label{Arrhenius}
 && c(x,\sigma) = c_1 (1 - \sigma(x)) + c_2\sigma(x)\EXP{-\beta U(x)}\COMMA\\
 && U(x) = K\sum_{y\in\Omega_x}\sigma(x+y) + h\COMMA
\end{eqnarray}
\begin{figure*}
\centerline{
    \subfigure[Coverage]{\label{phase1D}\includegraphics[width=0.4\textwidth]{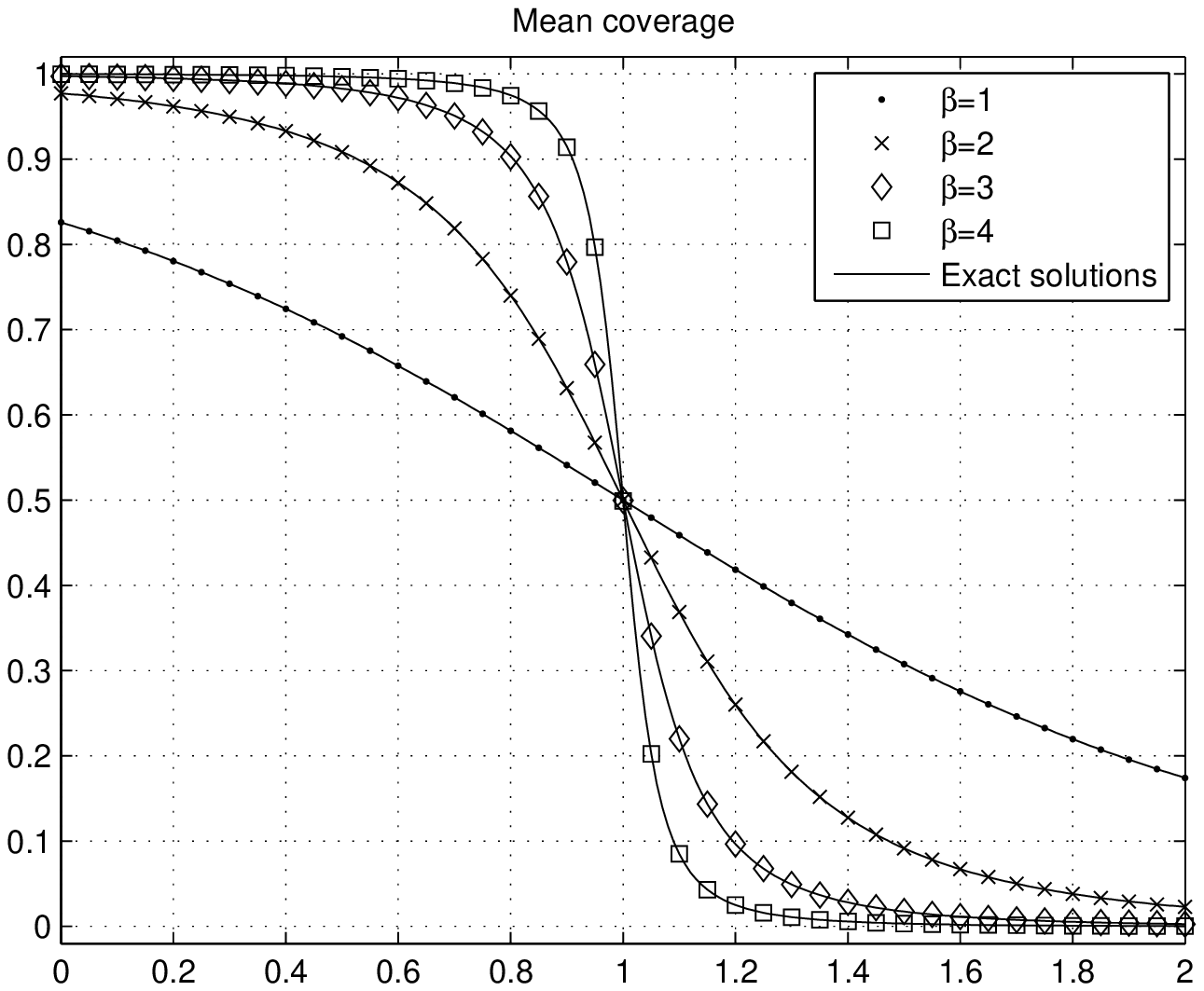}}
    \subfigure[Two-point spatial correlations]{\label{spcorr1D}\includegraphics[width=0.4\textwidth]{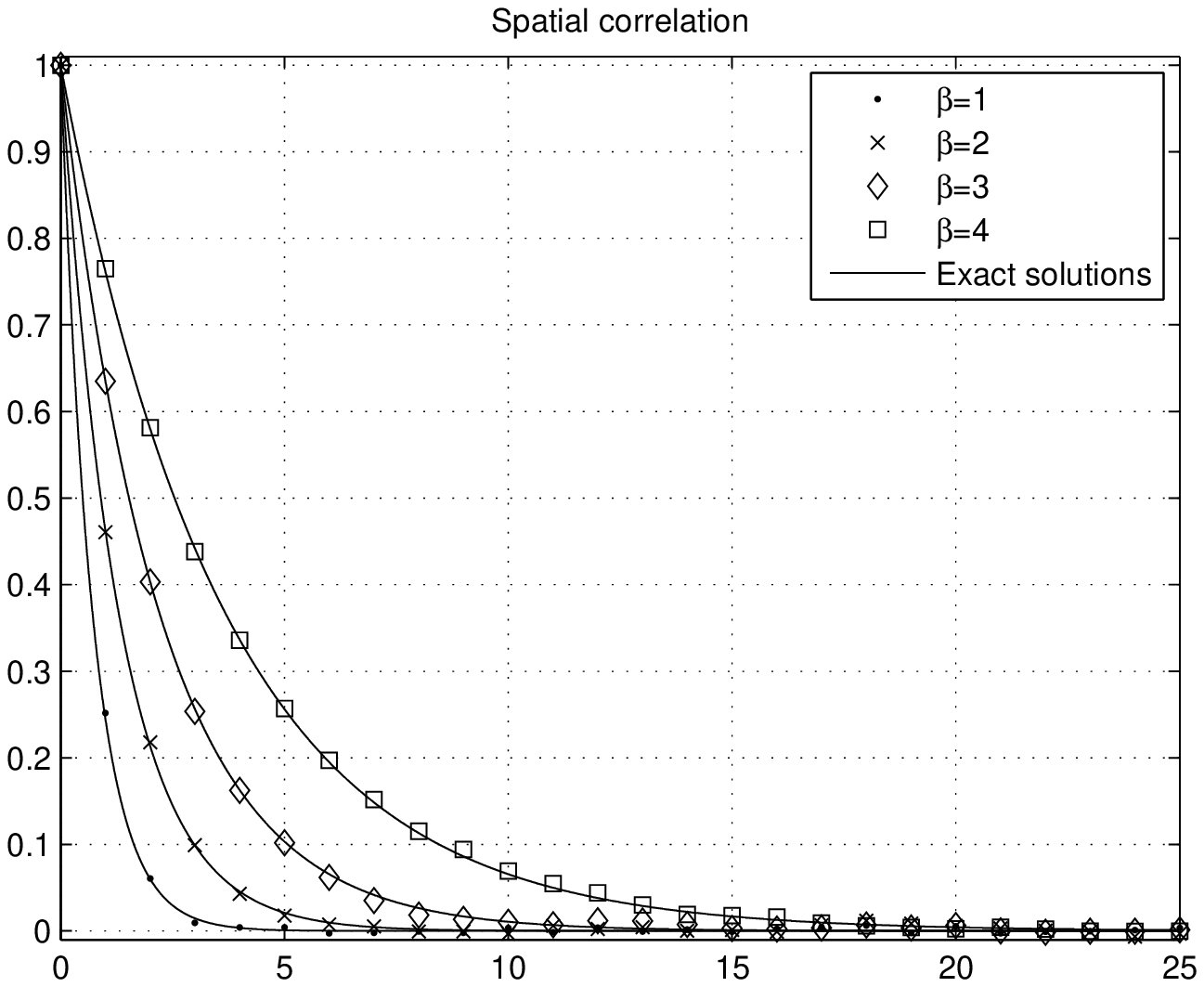}}
 }
\caption{(a)  Comparison of the exact solution \VIZ{exactcov1d} (solid line) for the total coverage $c_\beta(K,h)$, $K=1$, 
         with the mean coverage obtained in simulations
         on the one-dimensional lattice with $N=2^{15}$ and $\Delta t =1.0$. 
         (b) Two-point spatial correlation function
         estimated at $h=1$ on the same lattice and $\Delta t=1.0$ compared to the exact solution.}\label{phasediag1D} 
\end{figure*}
with $\beta$ is a given inverse temperature. The generator of \VIZ{Arrhenius} is a self-adjoint operator
on the space $\LTWO(\SIGMA,\mu_N)$ where $\mu_N(d\sigma) = Z^{-1}\EXP{-\beta H(\sigma)} \,d\sigma$ is the
canonical Gibbs measure of the system at the constant inverse temperature $\beta$. Consequently the dynamics
is reversible and the measure $\mu_t$ of the process $\PROC{S_t}$ converges to the Gibbs measure $\mu_N$ as
$t\to\infty$. Thus the dynamics \VIZ{Arrhenius} can be used for computing expected values 
$\EXPECT_{\mu_N}[f]$ by invoking ergodicity and averaging on a single trajectory
$$
  \EXPECT_{\mu_N}[f] \equiv \int_{\SIGMA} f(\sigma) \,\mu_N(d\sigma) = 
  \lim_{T\to\infty}\frac{1}{T}\int_0^T f(S_t) \,dt\PERIOD
$$
In the simulations we estimate two observables:
\begin{eqnarray*}
 &&\mbox{mean coverage:} \; \bar c_t = \frac{1}{|\LATT|}\EXPECT[\sum_{x\in\LATT}\sigma_t(x)]\COMMA\\
 &&\mbox{2-point correlation function:} \; \bar\lambda_t(x,y) = \EXPECT[\sigma_t(x)\sigma_t(x+y)]\PERIOD
\end{eqnarray*}
Due to translational invariance the function $\bar\lambda_k(x,y)$ depends on the distance
$|x-y|$ only. For exactly solvable one and two dimensional
Ising models we have explicit formulas which we summarize here for the spins in $\SPINSP=\{0,1\}$.

\smallskip
\noindent{\it 1D Ising model:} The one-dimensional Ising model does not exhibit a phase transition and thus presents
a simple benchmark for accuracy. Working with lattice gas models requires a simple transformation
of the well-known exact solution, \cite{Baxter}, 
which for the Hamiltonian of the system given on the periodic lattice
$$
H(\sigma) = -K\sum_{x=1}^N \sigma(x)\sigma(x+1) + h\sum_{x=1}^N\sigma(x)\COMMA
$$
\begin{figure}
\centerline{
\subfigure[Coverage]{\label{phase2D}\includegraphics[width=0.50\textwidth]{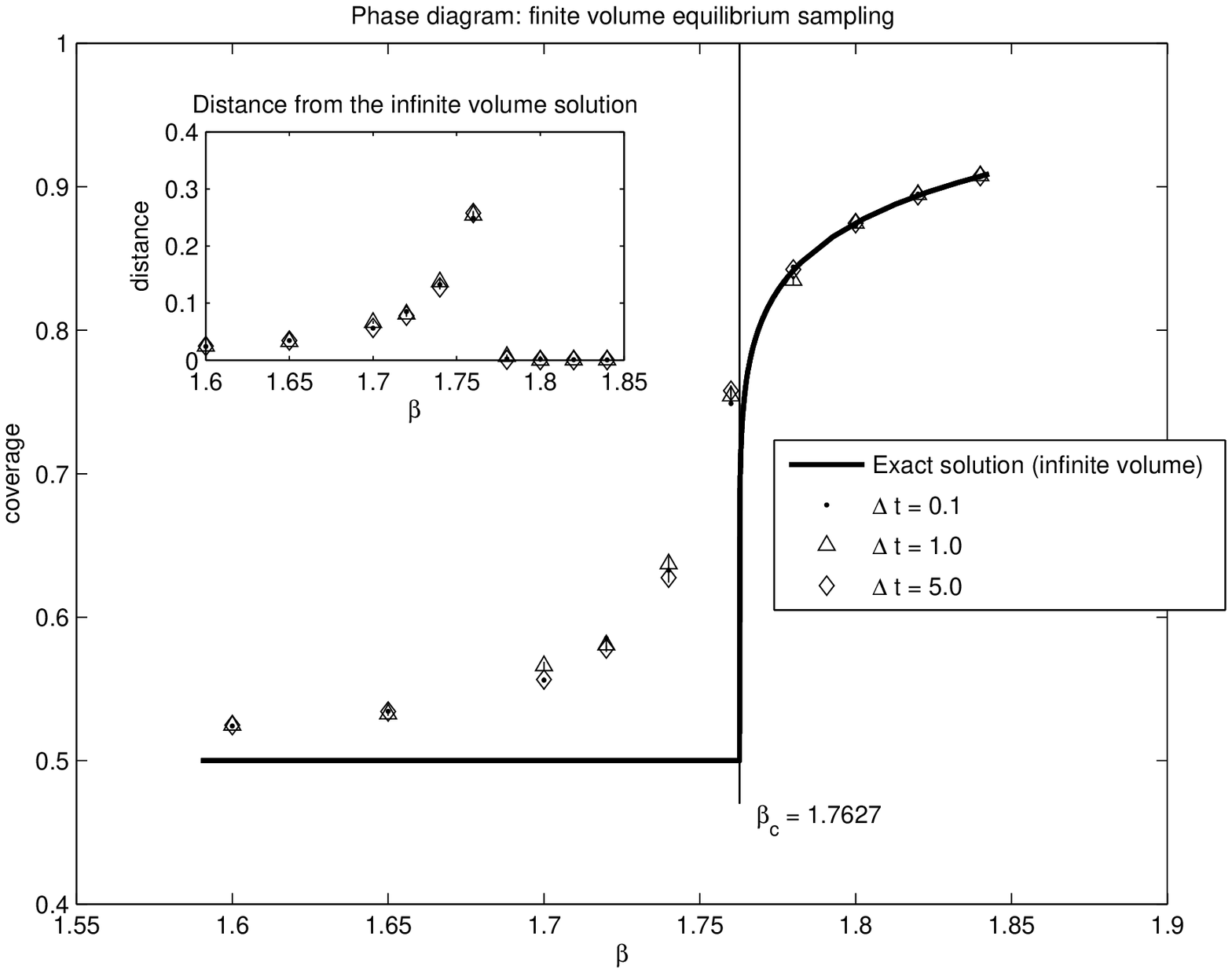}}
\subfigure[Two-point spatial correlations]{\label{2ptcorr}\includegraphics[width=0.50\textwidth]{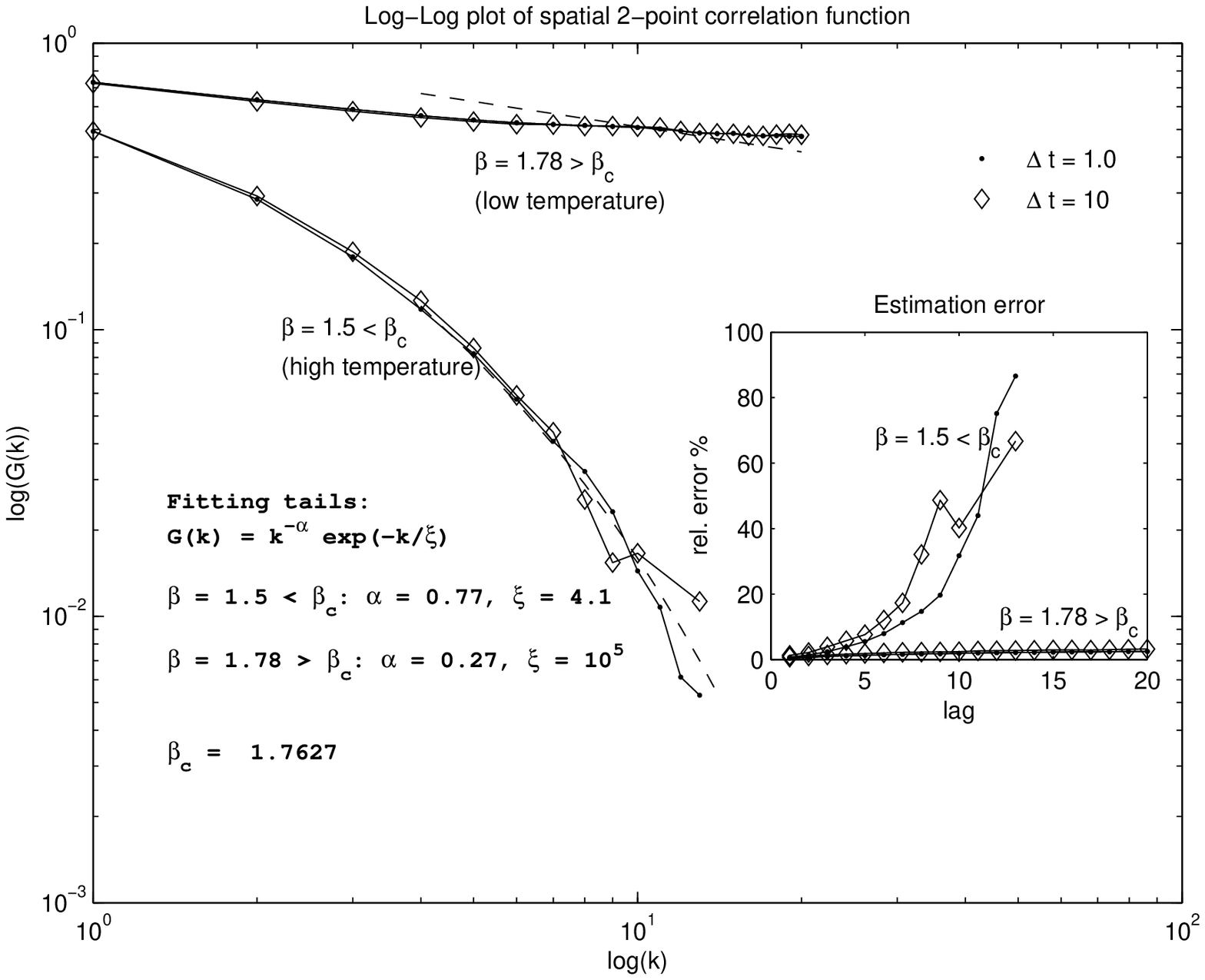}}
}
\caption{(a) Comparison of the exact solution \VIZ{exactcov2d} (solid line) for the total coverage $c_\beta(K,h)$, $h=2$, 
         with mean coverage obtained in simulations
         on the one-dimensional lattice with $N=128$ and various $\Delta t$'s. (b) Spatial two-point correlation function in the two-dimensional Ising model simulated
         on the lattice $N=512^{2}$ at a sub-critical temperature $\beta > \beta_c$ and supercritical
regime $\beta < \beta_c$. The simulation confirms the behavior obtained from the  infinite volume
exact solution: at high temperatures the decay is exponential while at temperatures below the
critical temperature the decay is algebraic. The dashed line represents the fitted function of
the form $k^{-\alpha}\EXP{-k/\xi}$.}
\end{figure}
yields the equilibrium mean coverage and the 2-point correlation function
\begin{align}
 \bar c(h,\beta) =& \frac{1}{2}\left( 1 +\frac{\SINH(h')}{\left(\SINH^2(h') + \EXP{-4K'}\right)^{1/2}}\right)\COMMA 
                     \\  \label{exactcov1d}
 \bar \lambda(x,y) =& \frac{1}{4}\left(1+\EXP{4 K'}\SINH^2(h')\right)\times\\
 &\left[\frac{\EXP{K'}\COSH(h') - \EXP{-K'}\left(1+\EXP{4K'}\SINH^2(h')\right)^{1/2}}%
             {\EXP{K'}\COSH(h') + \EXP{-K'}\left(1+\EXP{4K'}\SINH^2(h')\right)^{1/2}}\right]^{(x-y)}\COMMA\;
        y\geq x\COMMA \label{exactcorr1d}
\end{align}
where 
\begin{equation}
 K' = \frac{1}{4}\beta K \COMMA\;\;\;\mbox{and}\;\; h' = \frac{1}{2}\beta (h - K)\PERIOD 
\end{equation}
\begin{figure}
\centerline{
\subfigure{\label{pdf2d}\includegraphics[width=.5\textwidth]{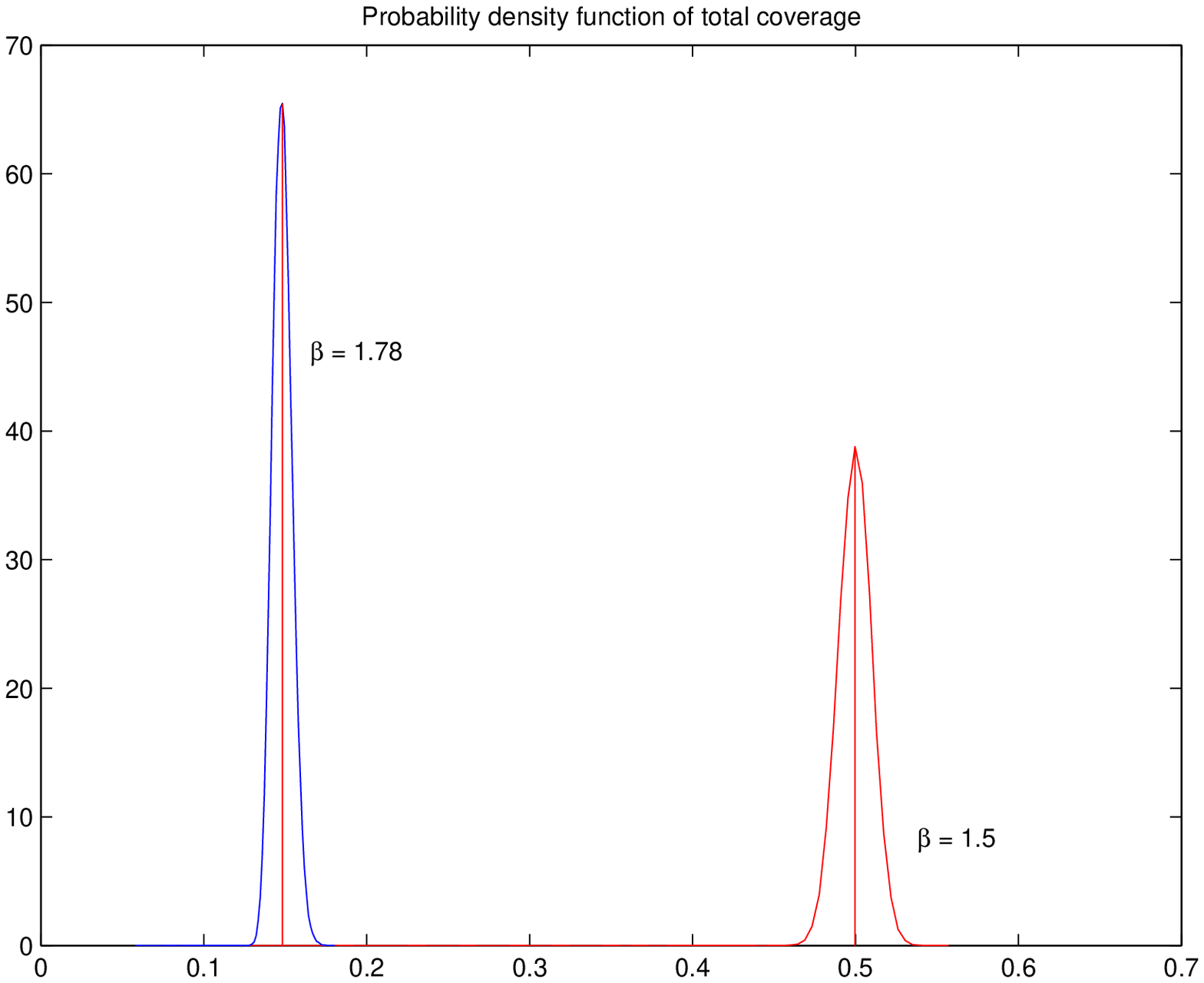}}
\subfigure{\label{dynamics2d}\includegraphics[width=0.5\textwidth]{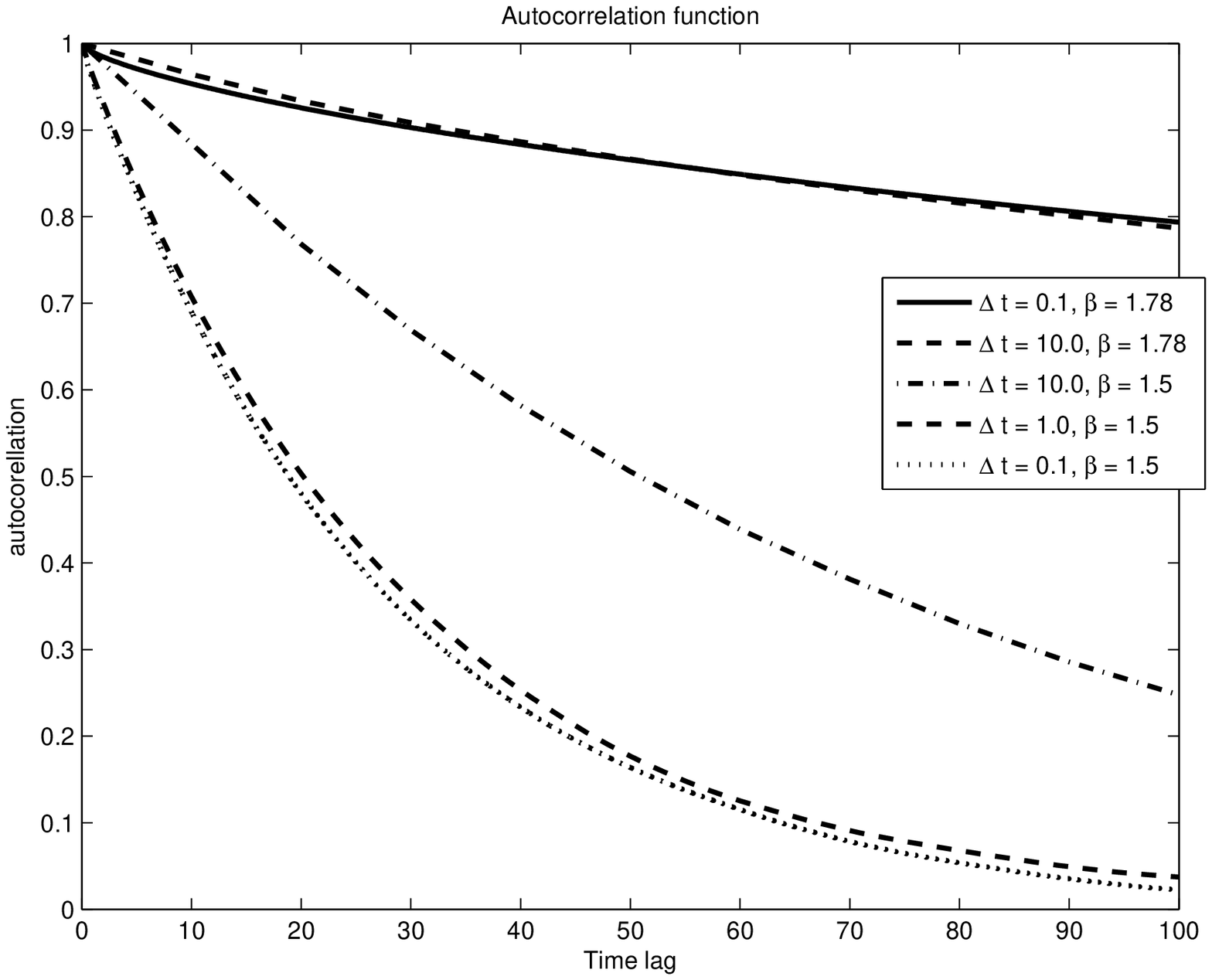}}
}
\caption{(a) Estimated equilibrium distributions of the coverage process at the two temperatures simulated in
         Fig.~\ref{2ptcorr}. (b) Autocorrelation functions for the coverage process in the two-dimensional Ising model
           simulated at $\beta = 1.5$ (high temperature above the critical temperature $\beta_c$ and
           at $\beta = 1.78 > \beta_c$ (low temperature), see parameters in  Fig.~\ref{2ptcorr}.
           }
\end{figure}
Since the one-dimensional Ising model does not exhibit a phase transition it allows us to assess the 
accuracy of the approximation for the phase diagram calculation. The phase diagram depicting dependence
of the coverage on the external field for different values of $\beta$ is shown in Figure~\ref{phase1D}.
In this simulation a rather conservative $\Delta t = 1.0 $ was chosen.
The statistical errors (confidence intervals) are below the resolution of the graph. 
As seen in the figure the isotherms for the average equilibrium coverage are thus obtained with 
a good accuracy. 
As a global observable the total coverage is less sensitive to statistical errors therefore we
also monitor the 2-point correlation function and its agreement with the exact solution \VIZ{exactcorr1d}. The results
for different values of $\beta$ in Figure~\ref{spcorr1D} demonstrate good accuracy.
\begin{figure}
\centerline{
    \subfigure[Sample path]{\label{path1D}\includegraphics[width=0.55\textwidth]{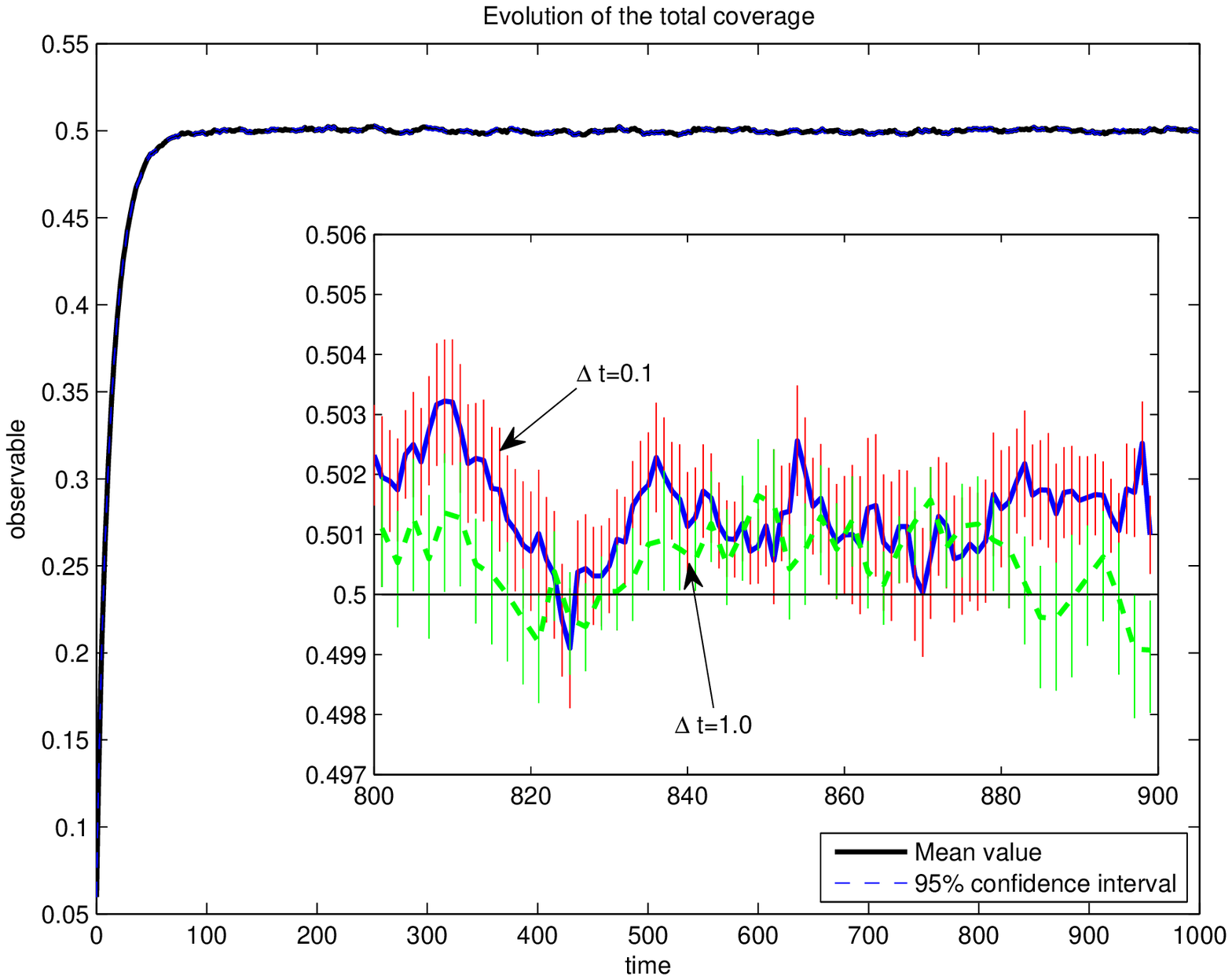}}
    \subfigure[Autocorrelation function]{\label{autocorr1D}\includegraphics[width=0.50\textwidth]{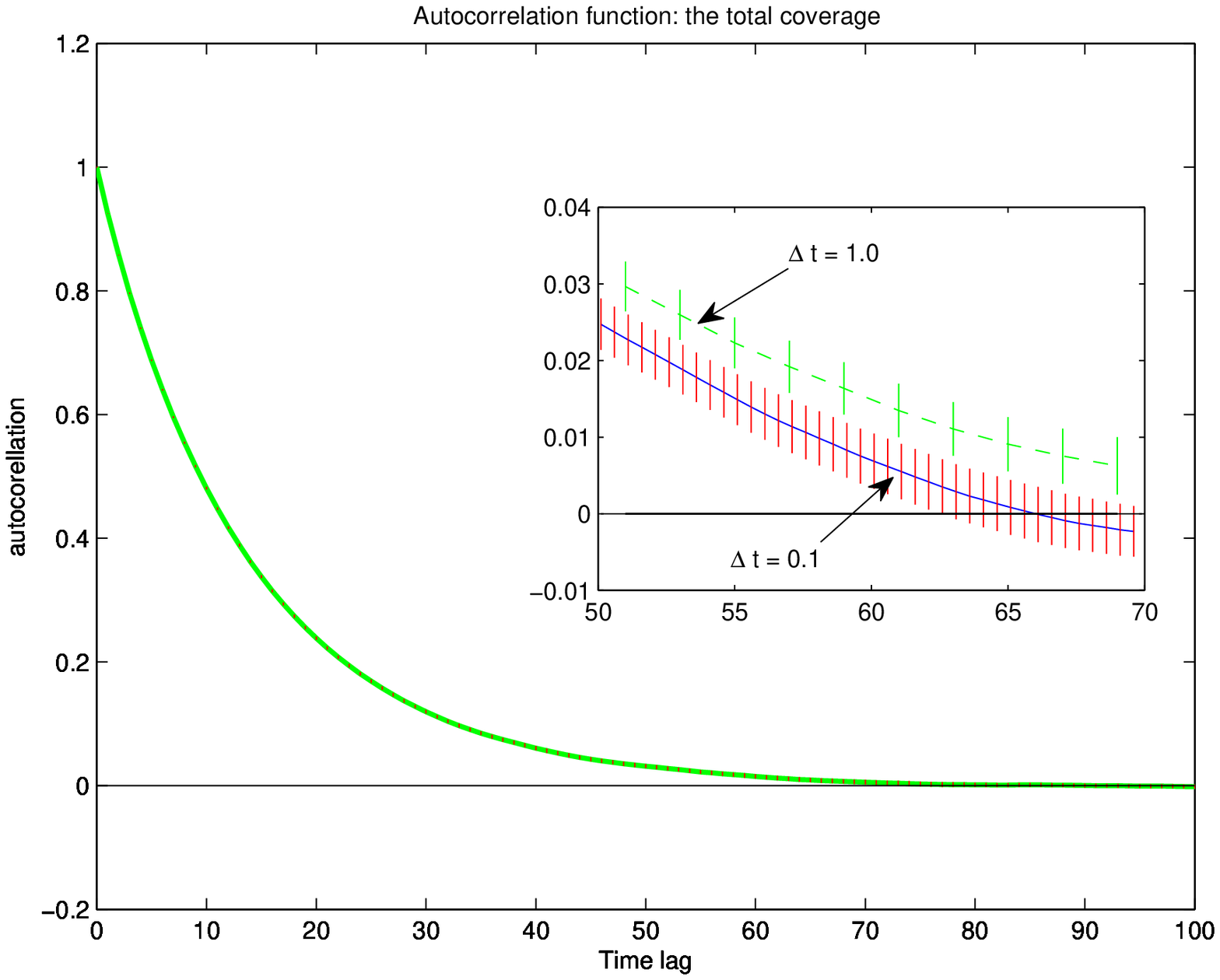}}
}
\caption{(a) A sample path of the total coverage process $\{S_t\}$ simulated at $\Delta t =1.0$ and $\Delta t = 0.1$ on the one-dimensional lattice with $N=2^{15}$ and $\Delta t =1.0$. 
         (b)  Autocorrelation function of the coverage process. The means were obtained from $M=1000$ independent
         realizations of the process at $\beta = 4$ and $h=1$. The inset shows error bars for the empirical mean estimator.}\label{samples1D} 
\end{figure}

\medskip
\noindent{\it 2D Ising model:} The phase transition that occurs in two-dimensional Ising model presents a more
challenging test case. However, the celebrated exact solution due to Onsager for spins $\SPINSP=\{-1,1\}$, 
\cite{Onsager44}, 
in the case with the zero external field and further
refinements yield closed formulas for the mean coverage and two point correlation functions. We restrict our
tests to the isotropic case, i.e., on the two-dimensional periodic lattice we have
the Hamiltonian
\begin{align*}
H(\sigma) = & - K \sum_{x=(x_1,x_2)\in\LATT} (\sigma(x_1,x_2) \sigma(x_1,x_2+1) + \\
            & \sigma(x_1,x_2)\sigma(x_1+1,x_2)) + h\sum_{x\in\LATT} \sigma(x)\PERIOD
\end{align*}
Transforming the exact solutions for the spins $\SPINSP=\{0,1\}$ we obtain the equivalent 
to the zero external field the value $h = 2K$ at which value
the critical inverse temperature solves $\SINH(\tfrac{1}{2}\beta_c K) = 1$. 
The exact solution for the mean coverage has the form
\begin{equation}\label{exactcov2d}
 \bar c(\beta) = 
  \begin{cases}   
     \frac{1}{2}\left(1 + \left[ 1 - (\SINH(\frac{1}{2}\beta K))^{-4}\right]^{1/8}\right) \COMMA & \beta > \beta_c\COMMA  \\
     \frac{1}{2}                                                                          \COMMA & \beta < \beta_c\PERIOD
  \end{cases}
\end{equation}
The exact solution for the 2-point correlation is available in \cite{WuMcCoy76}, however, we use only the asymptotics in
$|x-y|$, \cite{Baxter}. Introducing $\kappa = (\SINH(\tfrac{1}{2}\beta K))^{-2}$ we have  
\begin{equation}\label{exactcorr2d}
 \bar \lambda(x,y) = 
  \begin{cases}   
        (1 - \kappa^2)^{1/4} + \BIGO(\kappa^{|x-y|}) \COMMA & \beta > \beta_c\COMMA  \\
        \BIGO(\kappa^{-|x-y|/2})                     \COMMA & \beta < \beta_c\PERIOD
  \end{cases}
\end{equation}

The phase diagram is computed at $h=2$ which for $K=1$ corresponds to the regime when the 
second-order phase transition occurs at the critical temperature $\SINH(\tfrac{1}{2}K\beta_c)=1$. 
Sampling the coverage exhibits well-known difficulties close to the critical point $\beta_c$
which are not cured by the fractional step algorithm. Instead, we demonstrate in Figure~\ref{phase2D}
that for wide range of choices $\Delta t$ the phase diagram is constructed accurately for $\beta$
outside a neighborhood of $\beta_c$. Close to the critical point the algorithm provides approximations that
are in agreement with other Monte Carlo sampling approach. The finite-size effects are pronounced at
the neighborhood of the critical point due to algebraic decay of correlations. Thus it is not expected
that a good agreement with the infinite volume exact solution will be observed in the finite size 
simulations. Nonetheless, the presence of the second-order phase transition is indicated in the
computed phase diagram. Furthermore, the proposed algorithm provides an efficient implementation that allows
for simulations on large lattice. 
It is shown in Figure~\ref{2ptcorr} that
algebraic decay of the 2-point correlation function is well approximated in the low-temperature (sub-critical) 
regime, while at super-critical temperatures the exponential decay is observed.
Overall, we note that such long-time sampling of the simulated CTMC is a particularly challenging task since in principle, errors from any approximation may accumulate at long times and contaminate the simulation.

Studying approximation properties of the stochastic dynamics poses a more difficult task due to
the lack of an exact solution for the evolution of observables. Certain guidance can be obtained
from mean-field approximations, however, those do not give sufficiently good approximation for Ising
model in low dimensions. Therefore we compare the evolution of the coverage obtained from the traditional
SSA algorithm with approximations generated by the proposed fractional time step algorithm with
different choices $\Delta t$. In Figure~\ref{path1D} we compare the expected value and variance
of the total coverage process $C_t = \tfrac{1}{|\LATT|}\sum_{\LATT} S_t(x)$. 
Furthermore, it is also shown that the auto-correlation function for the process $C_t$ is well-approximated
and approximations converge as $\Delta t \to 0$, see Figures~\ref{autocorr1D} and \ref{dynamics2d}.

\subsection{Examples from Catalysis and  Reaction Engineering}

\begin{table}
\caption{An Event in $\Omega_x$, $x^{nn}\in\Omega_x$ is a randomly selected site from the nearest-neighbor set of $x$, and
$r_2(x) = \frac{1}{4}(1-\sigma(x)^2)\nu_0^x$, $r_3(x) = \frac{1}{8}\sigma(x)(1+\sigma(x))\nu_{-1}^x$, 
$r_4(x) = \frac{1}{8}\sigma(x)(\sigma(x)-1)\nu_{1}^x$, where $\nu_k^x$ is the number of nearest neighbors (n.n.)
of $x$ that are equal to $k$.
}\label{COrates}
\begin{tabular}{@{\vrule height 10.5pt depth4pt  width0pt}|l|lllll|}\hline
$\omega$& site         & $\sigma(x)$ & $\sigma^x$                                     & Rate $c(x,\omega;\sigma)$   &   Comment               \\ \hline
$1$     & vacant        & $0$  & $0\to  1$       & $k_1(1-(\sigma(x))^2)$ & $\mathrm{CO}$ adsorb  \\
$2$     & vacant        & $0$  & $0\to -1$       & $(1 - k_1)r_2(x)$      & $\mathrm{O}_2$ adsorb \\
        &               &      & $0 \to-1$, $x^{nn}$&                        &                \\
$3$     & $\mathrm{CO}$ & $1$  & $1\to 0$        & $k_2  r_3(x)$          & $\mathrm{CO}+\mathrm{O}$ and desorb \\
        &               &      & $-1 \to 0$, $x^{nn}$&                        &                \\
$4$     & $\mathrm{O}$  & $-1$ & $-1\to 0$       & $k_2  r_4(x)$          & $\mathrm{CO}+\mathrm{O}$ and desorb \\
        &               &      & $ 1\to 0$, $x^{nn}$&                        &                \\ \hline
\end{tabular}
\end{table}
In order to demonstrate the applicability of the proposed parallelization methodology in systems exhibiting complex spatio-temporal morphologies at mesoscopic length scales,e.g., islands, spirals, rings, etc., we implement a KMC algorithm arising in the modeling of chemical reaction dynamics  on a catalytic surface. Here we focus on CO oxidation, which is a prototypical example for molecular-level reaction-diffusion mechanisms between adsorbates on a surface.
We  note that molecular dynamics simulations have also been employed   to understand  micro-mechanisms  on surfaces such as  reaction paths \cite{Payne}. However, reaction kinetics for mesoscale adsorbate structures cannot be simulated by using molecular dynamics  because of spatio-temporal scale limitations of such methods, while  KMC methods, have the ability to simulate much larger scales \cite{Lukkien}.

In  KMC models for CO oxidation on a catalytic surface spatial resolution is a critical ingredient of the modeling since in-homogeneously
 adsorbed O and CO react on the catalytic surface only where the corresponding phases meet.
Sophisticated KMC models for CO oxidation on catalytic surfaces, where kinetic parameters are estimated by {\em ab initio} density functional theory (DFT), \cite{Kohn}, were recently developed in  \cite{Reuter1} and later in \cite{Ohta07}, \cite{evans09}. 
Such KMC models yield a remarkable agreement with experiments, 
see also  the review articles \cite{Metiu} and \cite{Christensen}. 

Next we demonstrate the performance of  parallel Fractional Step algorithms for KMC simulation  
to heterogeneous catalysis.  We implement a simplified CO oxidation model known as the  Ziff-Gulari-Barshad (ZGB) model, \cite{ZGB}, which 
was one of the first attempts towards a spatially distributed KMC modeling in reaction systems.
Although a simplified model compared to the {\em ab initio} KMC models described earlier, it incorporates the basic mechanisms 
for the dynamics
of adsorbate structures during CO oxidation on catalytic surfaces: single site updates (adsorption/desorption) and
multi-site updates (specifically, reactions with two sites being involved).
The spins take values $\sigma(x)=0$ denoting a vacant site $x\in\LATT$, $\sigma(x)=-1$ for a molecule 
$\mathrm{CO}$ at $x$,
and $\sigma(x)=1$ representing a $\mathrm{O}_2$ molecule. Depending on the local configurations 
of the nearest neighbors in $\Sigma_x = \{y\SEP |y-x|=1\}$ the events in Table~\ref{COrates} are executed.
The rates of individual events depend on 
the states in $\Omega_x$ which are enumerated by $\omega = \{1,2,3,4\}$ and are summarized in Table~\ref{COrates}.

\begin{figure}
\centerline{
\subfigure{\includegraphics[width=.4\textwidth]{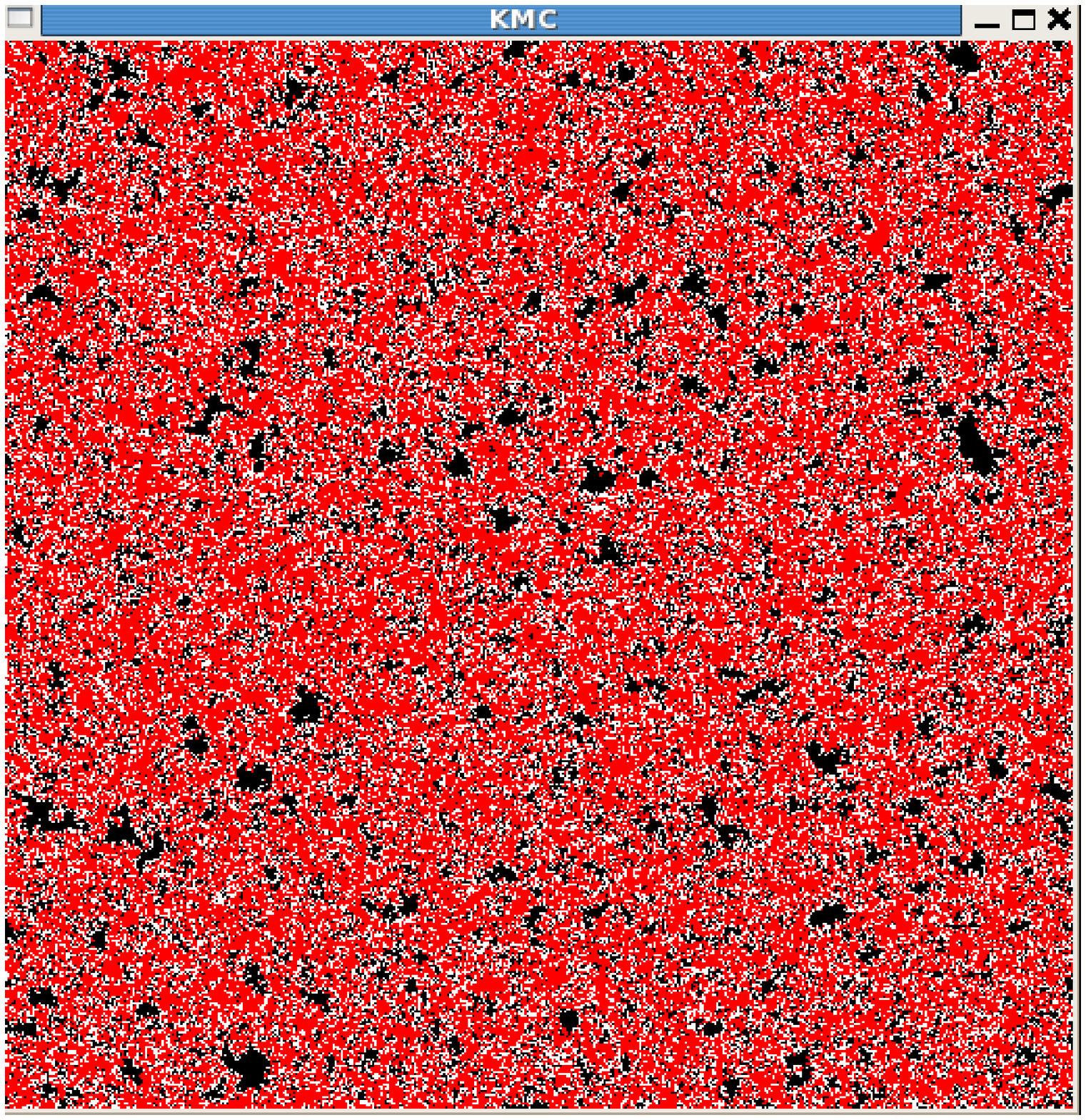}}
\subfigure{\includegraphics[width=.4\textwidth]{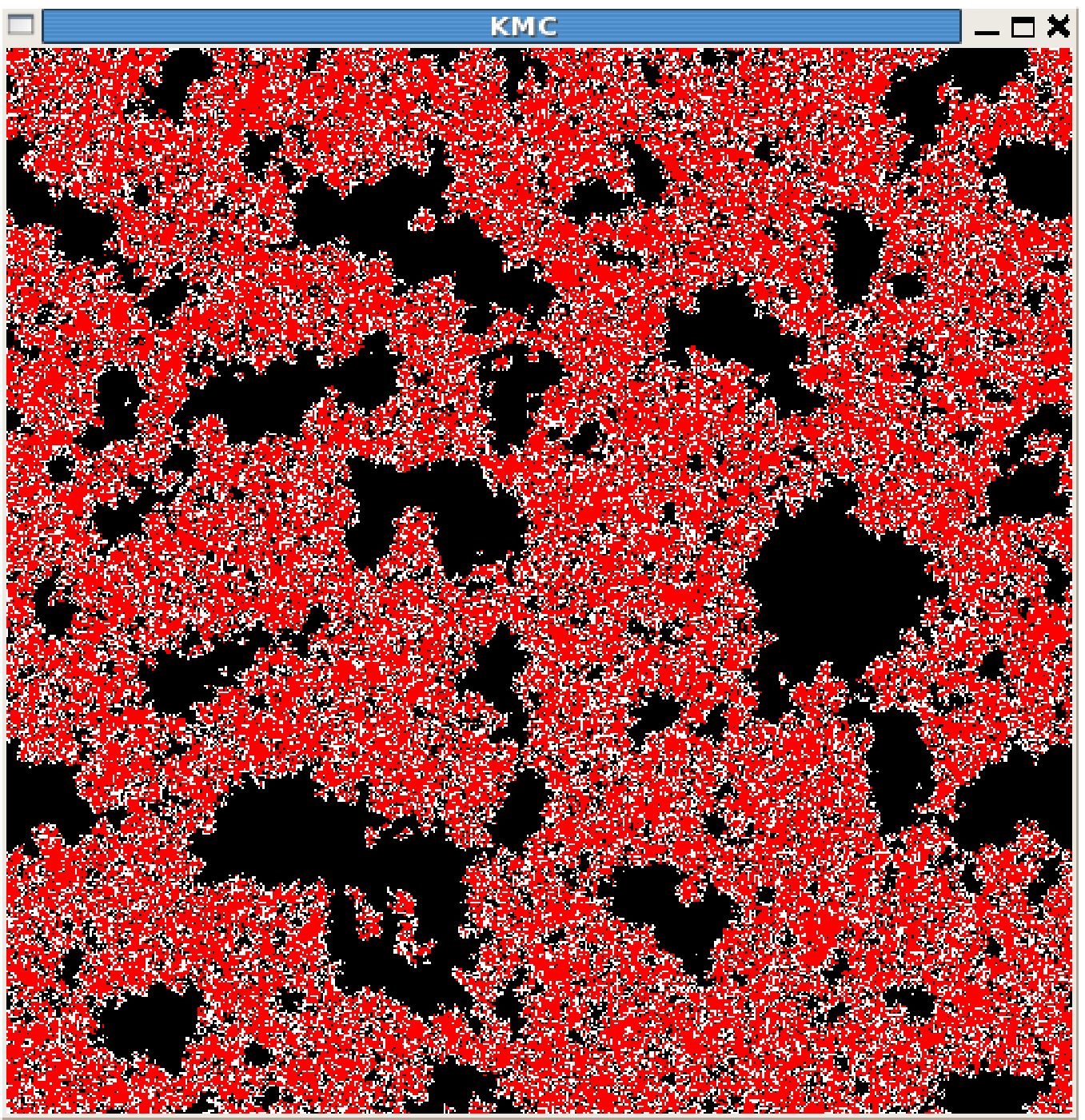}}
}
\centerline{\subfigure{\includegraphics[width=.4\textwidth]{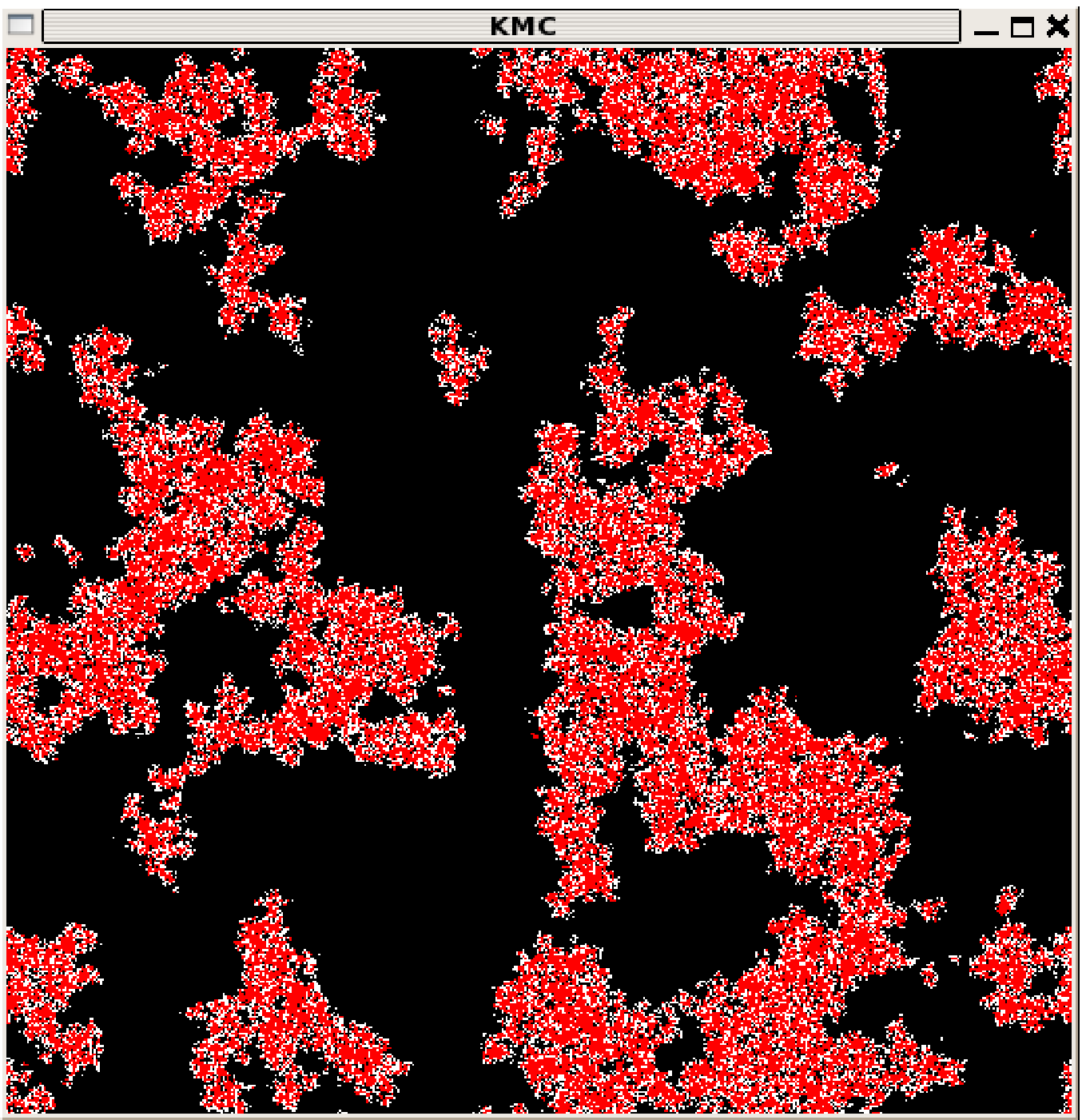}}
\subfigure{\includegraphics[width=.4\textwidth]{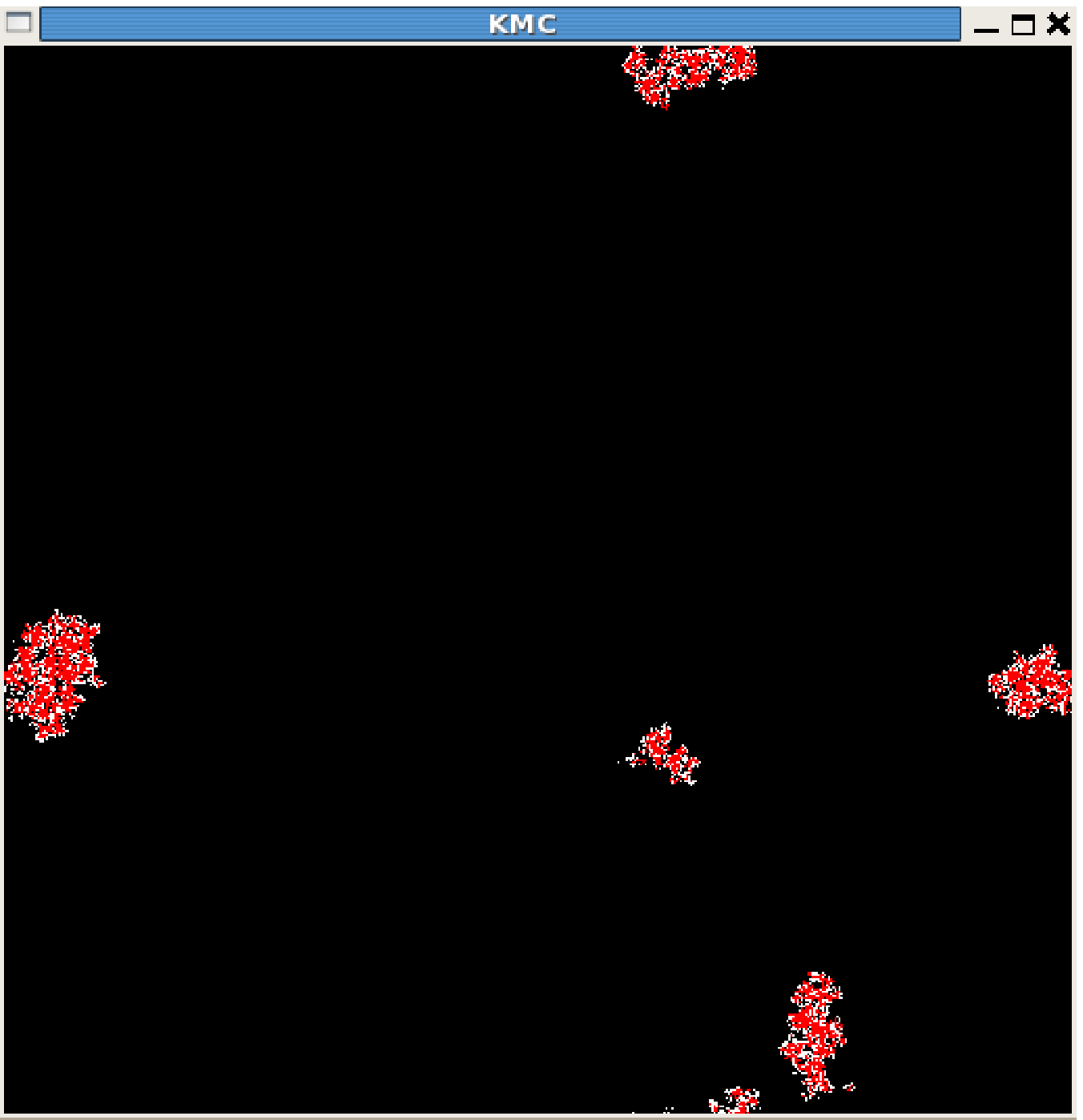}}
}
\caption{Snapshot at different simulation times for the $\mathrm{CO}$ oxidation process, 
         on a two-dimensional lattice $N=1024^2$.}\label{morphology}
\end{figure}
\begin{figure}
\centerline{\includegraphics[width=.8\textwidth]{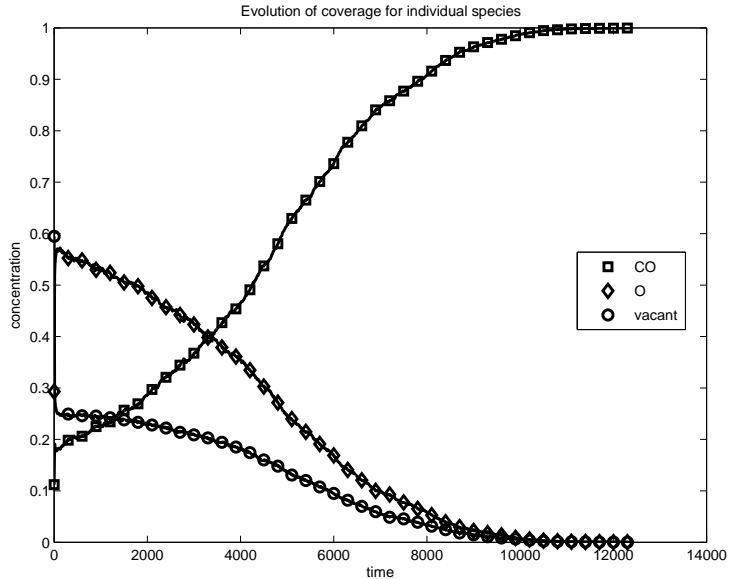}}
\caption{Evolution of the mean coverage for species in the oxidation process ($\mathrm{CO}$, $\mathrm{O}_2$,
and vacant sites).}\label{oxidation-coverage}
\end{figure}

The execution times for lattices of different sizes are compared in Figure~\ref{timeGPUs}, while a snapshot of the spatial
morphology is depicted in Figure~\ref{morphology}. 
Here we take as a reference the sequential KMC-BKL kernel. The same kernel is then used
for the implementation on GPUs where we compare times for different choices of $\Delta t$. 
We remark that the KMC kernel is not optimized by techniques such as the BKL algorithm, 
\cite{binder},
which is manifested in the scaling with respect to the size of the lattice $N$. However, the {\em same
kernel} is used in the fractional step algorithm thus we present fair comparisons between serial and parallel solvers,
noting that any optimized serial KMC algorithm can be used as a kernel in our Fractional Step framework.
It is worth noting that by partitioning 
of the problem into the subproblems the $\BIGO(N^2)$ complexity of the simple implementation for the SSA
algorithm is reduced, which is also demonstrated in Figure~\ref{timeGPUs} where the slope of lines
for simulations using GPUs suggest the  reduced complexity  of order $\BIGO(N)$.
Hence the proposed approach also offers a simple but
efficient implementation of KMC simulators.  

Finally, in our implementation (as well as in the original ZGB model) we did not implement the fast diffusion mechanism of $\mathrm{O}$ adsorbates on the surface, \cite{evans09}. 
However, the scheme \VIZ{strang} can allow us to easily implement within our parallelization framework
schemes with disparate time-scales which turn out to be important for the long-time adsorbate dynamics.

\section{Conclusions}

In this paper we proposed    a new framework for constructing parallel algorithms for lattice KMC simulations.
Our approach relies on a spatial decomposition of  the Markov generator underlying  the KMC  algorithm,  into  a hierarchy of operators corresponding  to  processors' structure  in the parallel architecture. Based on this operator 
decomposition, we  can formulate  Fractional Step Approximation schemes
by employing the Trotter product formula; these schemes  allow  us to run independently on each processor  a serial KMC simulation on each fractional time-step window.  Furthermore, the  schemes  incorporate  the  Communication Schedule between processors through 
the sequential application of the operators in the decomposition, as well as  the time step employed in the particular fractional step scheme.
Here we discussed  deterministic schedules resulting from Lie- and  Strang-type fractional step schemes, as well as random schedules derived 
by  the Random Trotter Theorem, \cite{Kurtz}.
We demonstrated   that the latter category includes the algorithm \cite{ShimAmar05b} as one particular example.

Some of the key features of the proposed  framework   and possible  future directions include:
The  hierarchical structure can be easily derived  and implemented for very general physiochemical processes modeled by lattice systems,  allowing    users to  input  as  the
KMC kernel their preferred serial algorithm. 
This flexibility and hierarchical structure allow for tailoring our framework to particular parallel architectures   with complex  memory and processor hierarchies, e.g., clusters of GPUs
communicating, for instance, through an MPI protocol,
 and
using the nested generator decomposition \VIZ{opdecomp2}.
Moreover, multi-scale Trotter algorithms for systems with  fast and slow processes are widely used  in Molecular Dynamics, e.g.,  \cite{Hairer}, and they can  be recast along  with the proposed  methods into a spatio-temporal hierarchy of operators that allow computational tasks to be hierarchically decomposed in space/time. 
The numerical consistency of the proposed algorithms is rigorously justified by Trotter Theorems, \cite{Trotter,Kurtz} showing the convergence 
of our approximating schemes to the original serial KMC algorithm.
Related numerical estimates are expected to provide insights on the design and the relative advantages
of various communication schedules  and architectures.
 We  discussed work load balancing between processors through  a re-balancing scheme based on probabilistic mass transport methods
that is particularly well-suited for the proposed fractional step KMC methods.
We carried out   detailed benchmarking  using analytically available exact solutions from statistical mechanics and applied  the method to simulate  complex spatially distributed molecular systems, such as reaction-diffusion processes on catalytic surfaces.
Finally, we studied the  performance and
scalability of our  algorithm \VIZ{opdecomp2} and the resulting code for different lattice sizes and
types of GPUs.

Concluding we note that there are some interesting conceptual analogies between the parallelization and coarse-graining algorithms of KMC such as
the Coarse-Grained Monte Carlo (CGMC) method
 e.g., \cite{KMV, AKPR}.  In both methods we decompose the particle system in components communicating minimally, e.g., \VIZ{opdecomp}, \VIZ{lie},
or trivially  as in coarse-graining methods,
thus, local information  is represented by  collective (coarse) variables, or computed on separate processors within a parallel architecture.
An early work towards parallelizing CGMC  \cite{KMV}
in problems with locally well-mixed  particle interactions is \cite{Taufer}, while further progress towards understanding and exploiting the 
analogies and the complementarity of CGMC and parallel KMC has the potential to give efficient KMC algorithms capable of simulating complex systems
at mesoscopic length scales.

\smallskip\noindent
{\bf Acknowledgments:}
The research of M.A.K. was partially supported by the National Science Foundation under the grant
NSF-DMS-071512, by  the Office of Advanced Scientific Computing Research,
U.S. Department of Energy under DE-SC0002339 and the European Commission FP7-REGPOT-2009-1 Award No 245749.
The research of P.P. was partially supported by the National Science Foundation under the grant
NSF-DMS-0813893 and by the Office of Advanced Scientific Computing Research,
U.S. Department of Energy under DE-SC0001340; the work was partly done at the Oak Ridge National Laboratory,
which is managed by UT-Battelle, LLC under Contract No. DE-AC05-00OR22725. 
The research of G.A. was partially supported
by the National Science Foundation under the grant
NSF-DMS-071512 and NSF-CMMI-0835582.
The research of M. T.  was partially supported by the AFOSR
STTR Program, the Army Research Office under the grant 54723-CS, and the  NVIDIA University Professor Partnership Program.
The research of L.X. was supported by the National Science Foundation under the grant NSF-DMS-0813893.


\begin{thebibliography}{39}
\expandafter\ifx\csname natexlab\endcsname\relax\def\natexlab#1{#1}\fi
\providecommand{\bibinfo}[2]{#2}
\ifx\xfnm\relax \def\xfnm[#1]{\unskip,\space#1}\fi
%
\bibitem[{Arampatzis et~al.(2011)Arampatzis, Katsoulakis, Plech{\'a}\v{c} and
  Rey-Bellet}]{AKPR11}
\bibinfo{author}{G.~Arampatzis}, \bibinfo{author}{M.A. Katsoulakis},
  \bibinfo{author}{P.~Plech{\'a}\v{c}}, \bibinfo{author}{L.~Rey-Bellet},
  \bibinfo{title}{Error analysis for parallel kinetic monte carlo algorithms:
  accuracy and processor communication}, \bibinfo{year}{2011}.
  \bibinfo{note}{Preprint}.
%
\bibitem[{Are et~al.(2008)Are, Katsoulakis, Plech{\'a}\v{c} and
  Rey-Bellet}]{AKPR}
\bibinfo{author}{S.~Are}, \bibinfo{author}{M.A. Katsoulakis},
  \bibinfo{author}{P.~Plech{\'a}\v{c}}, \bibinfo{author}{L.~Rey-Bellet},
  \bibinfo{title}{Multibody interactions in coarse-graining schemes for
  extended systems}, \bibinfo{journal}{SIAM J. Sci. Comput.}
  \bibinfo{volume}{31} (\bibinfo{year}{2008}) \bibinfo{pages}{987--1015}.
%
\bibitem[{Auerbach(2000)}]{Auerbach}
\bibinfo{author}{S.M. Auerbach}, \bibinfo{title}{Theory and simulation of jump
  dynamics, diffusion and phase equilibrium in nanopores.},
  \bibinfo{journal}{Int. Rev. Phys. Chem.} \bibinfo{volume}{19}
  (\bibinfo{year}{2000}).
%
\bibitem[{Battaile et~al.(2009)Battaile, Chandross, Holm, Thompson, Tikare,
  Wagner, Webb, Zhou, Cardona and Slepoy}]{SPPARKS}
\bibinfo{author}{S.P.C. Battaile}, \bibinfo{author}{M.~Chandross},
  \bibinfo{author}{L.~Holm}, \bibinfo{author}{A.~Thompson},
  \bibinfo{author}{V.~Tikare}, \bibinfo{author}{G.~Wagner},
  \bibinfo{author}{E.~Webb}, \bibinfo{author}{X.~Zhou}, \bibinfo{author}{C.G.
  Cardona}, \bibinfo{author}{A.~Slepoy}, \bibinfo{title}{Crossing the mesoscale
  no-man's land via parallel kinetic {M}onte {C}arlo}, \bibinfo{journal}{Sandia
  report}  (\bibinfo{year}{2009}).
%
\bibitem[{Baxter(1989)}]{Baxter}
\bibinfo{author}{R.J. Baxter}, \bibinfo{title}{Exactly Solved Models in
  Statistical Mechanics}, \bibinfo{publisher}{Academic Press},
  \bibinfo{edition}{3rd} edition, \bibinfo{year}{1989}.
%
\bibitem[{Bortz et~al.(1975)Bortz, Kalos and Lebowitz}]{BKL75}
\bibinfo{author}{A.B. Bortz}, \bibinfo{author}{M.H. Kalos},
  \bibinfo{author}{J.L. Lebowitz}, \bibinfo{title}{A new algorithm for {M}onte
  {C}arlo simulation of {I}sing spin systems}, \bibinfo{journal}{Journal of
  Computational Physics} \bibinfo{volume}{17} (\bibinfo{year}{1975})
  \bibinfo{pages}{10--18}.
%
\bibitem[{Chatterjee and Vlachos(2007)}]{Vlachos}
\bibinfo{author}{A.~Chatterjee}, \bibinfo{author}{D.~Vlachos},
  \bibinfo{title}{An overview of spatial microscopic and accelerated kinetic
  {M}onte {C}arlo methods}, \bibinfo{journal}{Journal of Computer-Aided
  Materials Design} \bibinfo{volume}{14} (\bibinfo{year}{2007})
  \bibinfo{pages}{253--308}. \bibinfo{note}{10.1007/s10820-006-9042-9}.
%
\bibitem[{Christensen and Norskov(2008)}]{Christensen}
\bibinfo{author}{C.H. Christensen}, \bibinfo{author}{J.K. Norskov},
  \bibinfo{title}{{A molecular view of heterogeneous catalysis}},
  \bibinfo{journal}{{Journal of Chemical Physics}} \bibinfo{volume}{{128}}
  (\bibinfo{year}{{2008}}).
%
\bibitem[{Eick et~al.(1993)Eick, Greenberg, Lubachevsky and
  Weiss}]{Lubachevsky93}
\bibinfo{author}{S.G. Eick}, \bibinfo{author}{A.G. Greenberg},
  \bibinfo{author}{B.D. Lubachevsky}, \bibinfo{author}{A.~Weiss},
  \bibinfo{title}{Synchronous relaxation for parallel simulations with
  applications to circuit-switched networks}, \bibinfo{journal}{ACM Trans.
  Model. Comput. Simul.} \bibinfo{volume}{3} (\bibinfo{year}{1993})
  \bibinfo{pages}{287--314}.
%
\bibitem[{Evans(1999)}]{lcevans}
\bibinfo{author}{L.C. Evans}, \bibinfo{title}{Partial differential equations
  and {M}onge-{K}antorovich mass transfer}, in: \bibinfo{booktitle}{Current
  developments in mathematics, 1997 ({C}ambridge, {MA})},
  \bibinfo{publisher}{Int. Press, Boston, MA}, \bibinfo{year}{1999}, pp.
  \bibinfo{pages}{65--126}.
%
\bibitem[{Gillespie(1976)}]{Gillespie76}
\bibinfo{author}{D.T. Gillespie}, \bibinfo{title}{A general method for
  numerically simulating the stochastic time evolution of coupled chemical
  reactions}, \bibinfo{journal}{Journal of Computational Physics}
  \bibinfo{volume}{22} (\bibinfo{year}{1976}) \bibinfo{pages}{403--434}.
%
\bibitem[{Hairer et~al.(2006)Hairer, Lubich and Wanner}]{Hairer}
\bibinfo{author}{E.~Hairer}, \bibinfo{author}{C.~Lubich},
  \bibinfo{author}{G.~Wanner}, \bibinfo{title}{Geometric numerical
  integration}, volume~\bibinfo{volume}{31} of
  \textit{\bibinfo{series}{Springer Series in Computational Mathematics}},
  \bibinfo{publisher}{Springer-Verlag}, \bibinfo{address}{Berlin},
  \bibinfo{edition}{second} edition, \bibinfo{year}{2006}.
  \bibinfo{note}{Structure-preserving algorithms for ordinary differential
  equations}.
%
\bibitem[{Heidelberger and Nicol(1993)}]{Nicol}
\bibinfo{author}{P.~Heidelberger}, \bibinfo{author}{D.M. Nicol},
  \bibinfo{title}{Conservative parallel simulation of continuous time {M}arkov
  chains using uniformization}, \bibinfo{journal}{IEEE Trans. Parallel Distrib.
  Syst.} \bibinfo{volume}{4} (\bibinfo{year}{1993}) \bibinfo{pages}{906--921}.
%
\bibitem[{Katsoulakis et~al.(2003)Katsoulakis, Majda and Vlachos}]{KMV}
\bibinfo{author}{M.~Katsoulakis}, \bibinfo{author}{A.~Majda},
  \bibinfo{author}{D.~Vlachos}, \bibinfo{title}{Coarse-grained stochastic
  processes for microscopic lattice systems}, \bibinfo{journal}{Proc. Natl.
  Acad. Sci} \bibinfo{volume}{100} (\bibinfo{year}{2003})
  \bibinfo{pages}{782--782}.
%
\bibitem[{Katsoulakis et~al.(2006)Katsoulakis, Plech{\'a}\v{c} and
  Sopasakis}]{KPS}
\bibinfo{author}{M.A. Katsoulakis}, \bibinfo{author}{P.~Plech{\'a}\v{c}},
  \bibinfo{author}{A.~Sopasakis}, \bibinfo{title}{Error analysis of
  coarse-graining for stochastic lattice dynamics}, \bibinfo{journal}{SIAM J.
  Numer. Anal.} \bibinfo{volume}{44} (\bibinfo{year}{2006})
  \bibinfo{pages}{2270--2296}.
%
\bibitem[{Kohn(1999)}]{Kohn}
\bibinfo{author}{W.~Kohn}, \bibinfo{title}{{Nobel {L}ecture: {E}lectronic
  structure of matter-wave functions and density functionals}},
  \bibinfo{journal}{{Reviews of Modern Physics}} \bibinfo{volume}{{71}}
  (\bibinfo{year}{{1999}}) \bibinfo{pages}{{1253--1266}}.
%
\bibitem[{Korniss et~al.(1999)Korniss, Novotny and Rikvold}]{Korniss99}
\bibinfo{author}{G.~Korniss}, \bibinfo{author}{M.A. Novotny},
  \bibinfo{author}{P.A. Rikvold}, \bibinfo{title}{Parallelization of a dynamic
  {M}onte {C}arlo algorithm: {A} partially rejection-free conservative
  approach}, \bibinfo{journal}{Journal of Computational Physics}
  \bibinfo{volume}{153} (\bibinfo{year}{1999}) \bibinfo{pages}{488--508}.
%
\bibitem[{Kurtz(1972)}]{Kurtz}
\bibinfo{author}{T.G. Kurtz}, \bibinfo{title}{A random {T}rotter product
  formula}, \bibinfo{journal}{Proc. Amer. Math. Soc.} \bibinfo{volume}{35}
  (\bibinfo{year}{1972}) \bibinfo{pages}{147--154}.
%
\bibitem[{Landau and Binder(2000)}]{binder}
\bibinfo{author}{D.P. Landau}, \bibinfo{author}{K.~Binder}, \bibinfo{title}{A
  guide to {M}onte {C}arlo simulations in statistical physics},
  \bibinfo{publisher}{Cambridge University Press},
  \bibinfo{address}{Cambridge}, \bibinfo{year}{2000}.
%
\bibitem[{Liggett(1985)}]{Liggett}
\bibinfo{author}{T.M. Liggett}, \bibinfo{title}{Interacting Particle Systems},
  volume \bibinfo{volume}{276} of \textit{\bibinfo{series}{Grundlehren der
  mathematischen {W}issenschaften}}, \bibinfo{publisher}{Springer-Verlag},
  \bibinfo{address}{New York, Berlin, Heidelberg, Tokyo}, \bibinfo{year}{1985}.
%
\bibitem[{Liu and Evans(2009)}]{evans09}
\bibinfo{author}{D.J. Liu}, \bibinfo{author}{J.W. Evans},
  \bibinfo{title}{{Atomistic and multiscale modeling of {CO}-oxidation on
  {P}d(100) and {R}h(100): {F}rom nanoscale fluctuations to mesoscale reaction
  fronts}}, \bibinfo{journal}{{Surf. Science}} \bibinfo{volume}{{603}}
  (\bibinfo{year}{{2009}}) \bibinfo{pages}{{1706--1716}}.
%
\bibitem[{Lubachevsky(1988)}]{Lubachevsky88}
\bibinfo{author}{B.D. Lubachevsky}, \bibinfo{title}{Efficient parallel
  simulations of dynamic {I}sing spin systems}, \bibinfo{journal}{J. Comput.
  Phys.} \bibinfo{volume}{75} (\bibinfo{year}{1988}) \bibinfo{pages}{103--122}.
%
\bibitem[{Lukkien et~al.(1998)Lukkien, Segers, Hilbers, Gelten and
  Jansen}]{Lukkien}
\bibinfo{author}{J.J. Lukkien}, \bibinfo{author}{J.P.L. Segers},
  \bibinfo{author}{P.A.J. Hilbers}, \bibinfo{author}{R.J. Gelten},
  \bibinfo{author}{A.P.J. Jansen}, \bibinfo{title}{{Efficient {M}onte {C}arlo
  methods for the simulation of catalytic surface reactions}},
  \bibinfo{journal}{{Physical Review E}} \bibinfo{volume}{{58}}
  (\bibinfo{year}{{1998}}) \bibinfo{pages}{{2598--2610}}.
%
\bibitem[{Mart{\'i}nez et~al.(2008)Mart{\'i}nez, Marian, Kalos and
  Perlado}]{Kalos}
\bibinfo{author}{E.~Mart{\'i}nez}, \bibinfo{author}{J.~Marian},
  \bibinfo{author}{M.H. Kalos}, \bibinfo{author}{J.M. Perlado},
  \bibinfo{title}{Synchronous parallel kinetic {M}onte {C}arlo for continuum
  diffusion-reaction systems}, \bibinfo{journal}{J. Comput. Phys.}
  \bibinfo{volume}{227} (\bibinfo{year}{2008}) \bibinfo{pages}{3804--3823}.
%
\bibitem[{Merrick and Fichthorn(2007)}]{MerickFichthorn07}
\bibinfo{author}{M.~Merrick}, \bibinfo{author}{K.A. Fichthorn},
  \bibinfo{title}{Synchronous relaxation algorithm for parallel kinetic {M}onte
  {C}arlo simulations of thin film growth}, \bibinfo{journal}{Phys. Rev. E}
  \bibinfo{volume}{75} (\bibinfo{year}{2007}) \bibinfo{pages}{011606}.
%
\bibitem[{Metiu(2008)}]{Metiu}
\bibinfo{author}{H.~Metiu}, \bibinfo{title}{{Preface to special topic: {A}
  survey of some new developments in heterogeneous catalysis}},
  \bibinfo{journal}{{Journal of Chemical Physics}} \bibinfo{volume}{{128}}
  (\bibinfo{year}{{2008}}).
%
\bibitem[{Nagasaka et~al.(2007)Nagasaka, Kondoh, Nakai and Ohta}]{Ohta07}
\bibinfo{author}{M.~Nagasaka}, \bibinfo{author}{H.~Kondoh},
  \bibinfo{author}{I.~Nakai}, \bibinfo{author}{T.~Ohta}, \bibinfo{title}{{CO}
  oxidation reaction on {Pt(111)} studied by the dynamic {M}onte {C}arlo method
  including lateral interactions of adsorbates}, \bibinfo{journal}{J. Chem.
  Phys.} \bibinfo{volume}{126} (\bibinfo{year}{2007})
  \bibinfo{pages}{044704--7}.
%
\bibitem[{Nandipati et~al.(2009)Nandipati, Shim, Amar, Karim, Kara, Rahman and
  Trushin}]{ShimAmar09}
\bibinfo{author}{G.~Nandipati}, \bibinfo{author}{Y.~Shim},
  \bibinfo{author}{J.G. Amar}, \bibinfo{author}{A.~Karim},
  \bibinfo{author}{A.~Kara}, \bibinfo{author}{T.S. Rahman},
  \bibinfo{author}{O.~Trushin}, \bibinfo{title}{Parallel kinetic {M}onte
  {C}arlo simulations of {Ag(111)} island coarsening using a large database},
  \bibinfo{journal}{Journal of Physics Condensed Matter} \bibinfo{volume}{21}
  (\bibinfo{year}{2009}) \bibinfo{pages}{084214}.
%
\bibitem[{Onsager(1944)}]{Onsager44}
\bibinfo{author}{L.~Onsager}, \bibinfo{title}{Crystal statistics. {I.} {A}
  two-dimensional model with an order-disorder transition},
  \bibinfo{journal}{Phys. Rev.} \bibinfo{volume}{65} (\bibinfo{year}{1944})
  \bibinfo{pages}{117--149}.
%
\bibitem[{Payne et~al.(1992)Payne, Teter, Allan, Arias and
  Joannopoulos}]{Payne}
\bibinfo{author}{M.C. Payne}, \bibinfo{author}{M.P. Teter},
  \bibinfo{author}{D.C. Allan}, \bibinfo{author}{T.A. Arias},
  \bibinfo{author}{J.D. Joannopoulos}, \bibinfo{title}{{Iterative minimization
  techniques for abinitio total-energy calculations - molecular-dynamics and
  conjugate gradients}}, \bibinfo{journal}{{Reviews of Modern Physics}}
  \bibinfo{volume}{{64}} (\bibinfo{year}{{1992}})
  \bibinfo{pages}{{1045--1097}}.
%
\bibitem[{Reuter et~al.(2004)Reuter, Frenkel and Scheffler}]{Reuter1}
\bibinfo{author}{K.~Reuter}, \bibinfo{author}{D.~Frenkel},
  \bibinfo{author}{M.~Scheffler}, \bibinfo{title}{{The steady state of
  heterogeneous catalysis, studied by first-principles statistical mechanics}},
  \bibinfo{journal}{{Physical Review Letters}} \bibinfo{volume}{{93}}
  (\bibinfo{year}{{2004}}).
%
\bibitem[{Sanders and Kandrot(2010)}]{CUDA}
\bibinfo{author}{J.~Sanders}, \bibinfo{author}{E.~Kandrot},
  \bibinfo{title}{{CUDA} by Example: {A}n Introduction to General-Purpose {GPU}
  Programming}, \bibinfo{publisher}{Addison-Wesley Professional},
  \bibinfo{address}{Cambridge}, \bibinfo{year}{2010}.
%
\bibitem[{Shim and Amar(2005{\natexlab{a}})}]{ShimAmar05}
\bibinfo{author}{Y.~Shim}, \bibinfo{author}{J.G. Amar},
  \bibinfo{title}{Rigorous synchronous relaxation algorithm for parallel
  kinetic {M}onte {C}arlo simulations of thin film growth},
  \bibinfo{journal}{Phys. Rev. B} \bibinfo{volume}{71}
  (\bibinfo{year}{2005}{\natexlab{a}}) \bibinfo{pages}{115436}.
%
\bibitem[{Shim and Amar(2005{\natexlab{b}})}]{ShimAmar05b}
\bibinfo{author}{Y.~Shim}, \bibinfo{author}{J.G. Amar},
  \bibinfo{title}{Semirigorous synchronous relaxation algorithm for parallel
  kinetic {M}onte {C}arlo simulations of thin film growth},
  \bibinfo{journal}{Phys. Rev. B} \bibinfo{volume}{71}
  (\bibinfo{year}{2005}{\natexlab{b}}) \bibinfo{pages}{125432}.
%
\bibitem[{Szabo and Fath(2007)}]{Szabo}
\bibinfo{author}{G.~Szabo}, \bibinfo{author}{G.~Fath},
  \bibinfo{title}{{Evolutionary games on graphs}}, \bibinfo{journal}{{Physics
  Reports}} \bibinfo{volume}{{446}} (\bibinfo{year}{{2007}})
  \bibinfo{pages}{{97--216}}.
%
\bibitem[{Trotter(1959)}]{Trotter}
\bibinfo{author}{H.F. Trotter}, \bibinfo{title}{On the product of semi-groups
  of operators}, \bibinfo{journal}{Proc. Amer. Math. Soc.} \bibinfo{volume}{10}
  (\bibinfo{year}{1959}) \bibinfo{pages}{545--551}.
%
\bibitem[{Wu et~al.(1976)Wu, McCoy, Tracy and Barouch}]{WuMcCoy76}
\bibinfo{author}{T.T. Wu}, \bibinfo{author}{B.M. McCoy}, \bibinfo{author}{C.A.
  Tracy}, \bibinfo{author}{E.~Barouch}, \bibinfo{title}{Spin-spin correlation
  functions for the two-dimensional {I}sing model: {E}xact theory in the
  scaling region}, \bibinfo{journal}{Phys. Rev. B} \bibinfo{volume}{13}
  (\bibinfo{year}{1976}) \bibinfo{pages}{316--374}.
%
\bibitem[{Xu et~al.(2010)Xu, Taufer, Collins and Vlachos}]{Taufer}
\bibinfo{author}{L.~Xu}, \bibinfo{author}{M.~Taufer},
  \bibinfo{author}{S.~Collins}, \bibinfo{author}{D.G.P. Vlachos},
  \bibinfo{title}{{Parallelization of Tau-Leap Coarse-Grained Monte Carlo
  Simulations on GPUs}}, in: \bibinfo{booktitle}{{Proceedings of the 2010
  IEEE/ACM International Parallel and Distributed Processing}}, {International
  Parallel and Distributed Processing Symposium (IPDPS)},
  \bibinfo{organization}{{}}, \bibinfo{publisher}{{in press}},
  \bibinfo{year}{{2010}}.
%
\bibitem[{Ziff et~al.(1986)Ziff, Gulari and Barshad}]{ZGB}
\bibinfo{author}{R.M. Ziff}, \bibinfo{author}{E.~Gulari},
  \bibinfo{author}{Y.~Barshad}, \bibinfo{title}{Kinetic phase transitions in an
  irreversible surface-reaction model}, \bibinfo{journal}{Phys. Rev. Lett.}
  \bibinfo{volume}{56} (\bibinfo{year}{1986}) \bibinfo{pages}{2553}.

\end{thebibliography}
\end{document}